\documentclass[final,1p]{elsarticle}
 \usepackage{epsfig}

\usepackage{amssymb}
\usepackage{color} %needed just for comments
\usepackage{placeins}

\journal{Journal of Computational Physics}

\def\vec#1{{\bf #1}}

\def\pfrac#1/#2.{\frac{\partial #1}{\partial #2}}
\def\ne{n_{\rm e}}
\def\np{n_{\rm p}}
\def\nn{n_{\rm n}}
\def\ve{\vec{v}_{\! \rm e}}
\def\vp{\vec{v}_{\! \rm p}}
\def\vn{\vec{v}_{\! \rm n}}
\def\De{D_{\rm e}}
\def\Dp{D_{\rm p}}
\def\Dn{D_{\rm n}}

\def\SeP{S_{\rm e}^{+}}
\def\SnP{S_{\rm n}^{+}}
\def\SpP{S_{\rm p}^{+}}
\def\SeM{S_{\rm e}^{-}}
\def\SnM{S_{\rm n}^{-}}
\def\SpM{S_{\rm p}^{-}}

\def\Rs{R_{\rm s}}
\def\xf{x_{\rm f}}
\def\xc{x_{\rm c}}
\def\epsi{\varepsilon_0}

\begin{document}

\begin{frontmatter}

%% Title, authors and addresses

%% use the tnoteref command within \title for footnotes;
%% use the tnotetext command for the associated footnote;
%% use the fnref command within \author or \address for footnotes;
%% use the fntext command for the associated footnote;
%% use the corref command within \author for corresponding author footnotes;
%% use the cortext command for the associated footnote;
%% use the ead command for the email address,
%% and the form \ead[url] for the home page:
%%
%% \title{Title\tnoteref{label1}}
%% \tnotetext[label1]{}
%% \author{Name\corref{cor1}\fnref{label2}}
%% \ead{email address}
%% \ead[url]{home page}
%% \fntext[label2]{}
%% \cortext[cor1]{}
%% \address{Address\fnref{label3}}
%% \fntext[label3]{}

\title{A new numerical strategy
       with space-time adaptivity and error control for multi-scale streamer
      discharge simulations\tnoteref{aknown}}

\tnotetext[aknown]{This research was 
supported by a fundamental project grant from ANR (French National Research Agency - ANR Blancs):
\emph{S\'echelles} (project leader S. Descombes), and by
a DIGITEO RTRA project: \emph{MUSE} (project leader M. Massot).
Authors express special thanks to Christian Tenaud (LIMSI-CNRS) 
for providing the basis of the multiresolution kernel of MR CHORUS, code 
developed for compressible Navier-Stokes equations (D\'eclaration d'Invention DI 03760-01).
}

%% use optional labels to link authors explicitly to addresses:
%% \author[label1,label2]{<author name>}
%% \address[label1]{<address>}
%% \address[label2]{<address>}

 \author[em2c,ecp]{Max Duarte\corref{cor1}\fnref{max_s}}
 \ead{max.duarte@em2c.ecp.fr}
 \fntext[max_s]{Ph.D. grant
from Mathematics (INSMI) and Engineering (INSIS) Institutes of CNRS 
and supported by INCA project (National Initiative for Advanced Combustion - CNRS - ONERA - SAFRAN).}

\author[em2c,ecp,brno]{Zden\v{e}k Bonaventura\fnref{zdenek_s}}
 \ead{zbona@physics.muni.cz}
\fntext[zdenek_s]{Support of Ecole Centrale Paris is gratefully acknowledged
for several month stay of Z. Bonaventura at Laboratory EM2C
as visiting Professor. Z. Bonaventura is also grateful to the 
Ministry of Education, Youth and Sports of the Czech Republic under project 
CZ.1.05/2.1.00/03.0086 and project MSM 0021622411.}

\author[em2c,ecp]{Marc Massot}
 \ead{marc.massot@em2c.ecp.fr}

\author[em2c,ecp]{Anne Bourdon}
 \ead{anne.bourdon@em2c.ecp.fr}

\author[nice]{St\'ephane Descombes}
\ead{sdescomb@unice.fr}

 \author[icj]{Thierry Dumont}
 \ead{tdumont@math.univ-lyon1.fr}

\cortext[cor1]{Corresponding author}

\address[em2c]{CNRS, UPR 288 ``Laboratoire d'Energ\'etique Mol\'eculaire et Macroscopique, Combustion'',
Grande voie des vignes, 92295
Ch\^{a}tenay-Malabry, France}
\address[ecp]{Ecole Centrale Paris, Grande Voie des Vignes, 92295 Ch\^{a}tenay-Malabry Cedex,
 France}
\address[brno]{Department of Physical Electronics, Faculty of Science, Masaryk
University, Kotl\'a\v{r}sk\'a 2, 611 37 Brno, Czech Republic}
 \address[nice]{Laboratoire J. A. Dieudonn\'e -
 UMR CNRS 6621,
 Universit\'e Nice Sophia Antipolis,
 Parc Valrose,
 06108 Nice Cedex 02,
 France}
 \address[icj]{Universit\'e de Lyon,
Universit\'e Lyon 1,
INSA de Lyon, 
Ecole Centrale de Lyon,
Institut Camille Jordan - UMR CNRS 5208,
43 Boulevard du 11 novembre 1918,
69622 Villeurbanne Cedex, France}

\begin{abstract}
%% Text of abstract
This paper presents a new resolution strategy for multi-scale streamer discharge simulations 
based on a second order time adaptive integration and space adaptive multiresolution. 
A classical fluid model is used to describe plasma discharges, 
considering drift-diffusion equations and the computation of electric field. 
The proposed numerical method provides a time-space accuracy control of the solution,
and thus, an effective accurate resolution independent of the fastest physical time scale.
An important improvement of the computational efficiency is achieved whenever the required 
time steps go beyond standard stability constraints
associated with mesh size or source time scales
for the resolution of the drift-diffusion equations,
whereas the
stability constraint related to the dielectric relaxation time scale
is respected but with a second order precision.
Numerical illustrations show
that the strategy can be efficiently applied to simulate  
the propagation of highly nonlinear ionizing waves as streamer discharges,
as well as highly multi-scale
nanosecond repetitively pulsed discharges,
describing consistently
a broad spectrum of space and time scales 
as well as different physical scenarios
for consecutive discharge/post-discharge phases,
out of reach of standard non-adaptive methods. 
\end{abstract}
\begin{keyword}
multi-scale discharge \sep time adaptive integration \sep space adaptive multiresolution \sep error control
%% keywords here, in the form: keyword \sep keyword

%% MSC codes here, in the form: \MSC code \sep code
%% or \MSC[2008] code \sep code (2000 is the default)
\MSC 65M08 \sep 65M50 \sep 65Z05 \sep 65G20
\end{keyword}

\end{frontmatter}

%%
%% Start line numbering here if you want
%%
% \linenumbers
%%%%%%%%%%%%%%%%%%%%%%%%%%%%%%%%%%%%%%%%%%%%%%%%%%%%%%%%%%%%%%%%%%%%%%%%%
%% main text
\section{Introduction}
\label{SecIntro}
In recent years, plasma discharges at atmospheric pressure have been studied for an increasing list of applications 
such as chemical and biological decontamination, aerodynamic flow control and combustion \cite{vanVeldhuizen:2000,Fridman:2005}.
In all these physical configurations, the discharges take usually the form of thin plasma filaments driven by highly 
nonlinear ionizing waves, also called streamers. 
These ionizing waves occur as a consequence of the high electric field induced by the fast variations of the net 
charge density ahead of an electron avalanche with large amplification.
The streamer discharge dynamics 
are mainly governed by the Courant, the effective ionization and the dielectric relaxation times scales \cite{Vitello:1994},
which are usually  of the order of $10^{-14}-10^{-12}$s,
whereas the
typical time scale of the discharge propagation
in centimeter gaps, is about a few tens of nanoseconds. 
On the other hand,
a large variation of space scales needs also to be taken into account, 
since the Debye length at atmospheric pressure can be as small as a few micrometers, 
while the inter-electrode gaps, where discharges propagate, are usually of the order of a few centimeters.
As a result, 
the detailed physics of the discharges reveals an important time-space 
multi-scale character 
\cite{Unfer:2010,Ebert_nonlinearity:2011}.

More complex applications include plasma assisted combustion or flow control, 
for which the enhancement of the gas flow chemistry or momentum transfer 
during typical time scales of the flow of $10^{-4}-10^{-3}$s, 
is due to consecutive discharges generated by high frequency (in the kHz range) sinusoidal or 
pulsed applied voltages \cite{Pilla:2006,opaits:2008}. 
Therefore, during the post-discharge phases of the order of tens of microseconds, not only 
the time scales are very different from those of the discharge phases of a few tens of nanoseconds,
but a completely different physics is taking place.
Then, to the rapid multi-scale configuration during discharges,
we have to add other rather slower multi-scale phenomena 
in the post-discharge,  
such as recombination of charged species, heavy-species chemistry, diffusion, gas heating and convection.
Therefore,
it is 
very challenging 
to accurately simulate the physics of plasma/flow interactions due to the synergy effects between 
the consecutive discharge/post-discharge phases.

In most numerical models of streamer discharges, 
the motion of electrons and ions is governed by drift-diffusion equations coupled with Poisson's equation. 
Early simulation 
studies were limited to simplified situations where the streamer is
considered as a cylinder of constant radius \cite{Davies:1964,Davies:1971,Abbas:1980,Morrow:1985},
in which the charged particle densities are 
assumed to be constant along the radial extension of the streamer:
the {\it 1.5D model} approach. 
In this model, the spatio-temporal 
evolution of the charged particle densities is solved only along one
spatial dimension in the direction of propagation,
whereas the
electric field  is calculated in two dimensions
using the so-called {\it disc method},
based on a direct integration of analytical results.
A 2D model for the electric field is indeed essential
to properly calculate the electric field enhancement by the space
charge in the streamer head.
After the first 2D streamer simulations using 
the Poisson's equation resolution were performed
\cite{Dhali:1987}, 
many studies have been carried out in 2D
\cite{Vitello:1994,Babaeva:1997,
Kulikovsky:2000,Pancheshnyi:2001,Arrayas:2002,Celestin:2009,Bourdon:2010} 
and 3D \cite{Nikandrov:2008,Pancheshnyi:2008,Luque_PRL:2008,Papageorgiou:2010}.

Being aware of the complexity of fully coupled resolutions of these
modeling equations,
a decoupling strategy is usually adopted,
which considers an independent and successive numerical resolution of 
Poisson's equation with a fixed charge distribution,
and of the drift-diffusion equations with a fixed electric field during each decoupling time step.
These computations might be performed explicitly in time 
with standard first or even second order schemes \cite{Montijn:2006,Bourdon:2007}.
In these cases,
the time steps are usually limited
for the sake of stability by the various characteristic times scales
(Courant, ionization, dielectric relaxation),
whereas the accuracy of simulations is assumed to be given by the resolution
of the fastest physical time scale.
In order to somehow overcome the dielectric relaxation limitation,
some semi-implicit approaches were developed \cite{Ventzek1993,Colella1999,Hagelaar:2001},
based on a predictive approximation of the space charge ahead in time 
during the electric field computation,
even though the other time scale constraints remain.
This gain of stability allows important improvements in terms of computational efficiency
but the accuracy of simulations becomes rather difficult to quantify.

Another performing technique to improve the efficiency of
simulations considers an asynchronous explicit time integration of the
drift-diffusion equations with self-adaptive local time-stepping,
for which the local time steps are based either on local dynamic increments of the solution 
\cite{Karimabadi05,Omelchenko06}
or on local Courant conditions \cite{Unfer07}.
These techniques are the subject of several studies
\cite{Coquel2008,Domingues2008,coquel:540}
and are mainly conceived to avoid expensive computations
whenever the whole system is unnecessarily advanced in time with a global time step prescribed
by the fastest scale.
Even though these methods yield efficient strategies,
specially in terms of CPU time savings, with
stable and flux-conserving time integrations,
it is rather difficult to conduct 
an accuracy control on the resolution of the
time dependent equations
or on their coupling with the electric field resolution for plasma models.

In this work, a numerical study is conducted in order to build
a second order explicit in time decoupling scheme for the resolution of
the electric field and the electron and ion densities.
A lower order and embedded method is taken into account to dynamically
compute the decoupling time steps that guarantee an accurate description
with error control of the global physical coupling.
At this stage, the only limiting time scale is the dielectric relaxation
characteristic time 
for stability reasons.
In a second level,
the drift-diffusion equations are solved using a Strang second order operator splitting
scheme in order to guarantee the global order of the strategy
\cite{article_mr,article_avc}.
This time integration scheme considers
high order dedicated methods during each
splitting time step, which is dynamically adapted by an error control procedure \cite{article_pas}.
In this way,
even though there is a global advance in time
given by the splitting time step,
the latter is
determined by the
desired accuracy of the global physics, 
which is not
necessarily related to 
the stability constraints
associated with the mesh size or the fastest source time scales
as demonstrated in \cite{article_fvca}.
As a consequence,
this technique provides an error control procedure
and stands as an alternative way to local stepping schemes
to overcome time step limitations related to
the reaction, diffusion and convection phenomena.

Both the electric field and density resolutions are performed
on an adapted mesh obtained by
a spatial multiresolution method,
based on Harten's pioneering work \cite{Harten95}
and further developed in \cite{Cohen03},
taking into account the
spatial multi-scale features of
these phenomena with
steep spatial gradients.
In particular, some grid adaptation techniques for 2D structured meshes were already used 
\cite{Montijn:2006,Unfer:2010,Pancheshnyi:2008} and extensions to 3D have been 
also proposed \cite{Pancheshnyi:2008,Nikandrov:2008}
for streamer simulations.
However,
one of the main advantages of the multiresolution approach is
that it is based on a wavelet representation technique and an error of the spatial approximation can be then
mathematically estimated.
Consequently,
an effective error control is achieved
for both the time and space resolution of the multi-scale phenomena under study.

The performance of the method is first evaluated for a 
propagating streamer problem with the multi-scale features previously discussed, 
for which the various simulation parameters are 
studied.
Once the physical configuration is settled, a 
1.5D streamer model is adopted in order to 
obtain 
an electric field resolution strategy
based on direct computations and
derived from analytical expressions,
suitable for adapted
finite volume discretizations \cite{Bessieres:2007}.
In a second step,
a more complex physical configuration is considered 
for the simulation of repetitively pulsed discharges, for which
a time-space adaptive method is required to efficiently
overcome some highly multi-scale
features in order to fully describe the various physical phenomena.
In this work, only a 1.5D model is considered 
but extensions to higher dimensions is straightforward for instance 
with a Poisson's 
equation solver for adapted grids as it has been implemented in 
\cite{Montijn:2006,Unfer:2010,Pancheshnyi:2008}.
However, in this paper we focus on the development 
and validation of
new numerical methods
for the resolution of the drift-diffusion equations
and its coupling with the electric field computation,
which are independent of the dimension of the problem.
Numerical illustrations of multidimensional problems
with the same time-space adaptive strategy with error control
will be the subject of future work.

The paper is organized as follows:
in Section~\ref{sec:model}, we present the physical configuration and the modeling equations.
The numerical strategy is presented in Section~\ref{sec:num_str}, in which the
second order adaptive time integration technique is detailed
along with the resolution of drift-diffusion equations and the electric field,
as well as the spatial multiresolution adaptive procedure.
Numerical illustrations are summarized in Section~\ref{sec:num_res}
for two configurations given by single propagating and multi-pulsed discharges.
We end with some concluding remarks and prospects on future developments
and applications.

%%%%%%%%%%%%%%%%%%%%%%%%%%%%%%%%%%%%%%%%%%%%%%%%%%%%%%%%%%%%%%%%%%%%%%%%%%%%%%%%%%%%%%%%%%%%
\section{Model formulation}\label{sec:model}
\begin{figure}[!htb]
 \begin{center}
   \includegraphics[width=0.5\hsize]{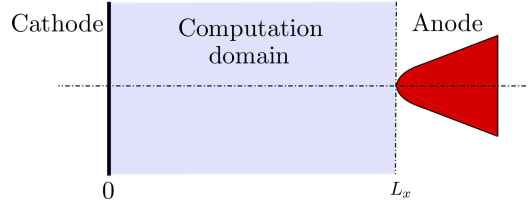}
   \caption{Computational domain for the studied point-to-plane geometry.
 \label{fdomain}}
 \end{center}
\end{figure}

In this work, we consider positive streamer discharges in air at atmospheric pressure in 
a point-to-plane geometry, as shown in Figure~\ref{fdomain}.
The tip of the anode is placed 
$1\,$cm from the planar cathode and the radius of curvature
of the anode is $324\,\mu$m. 
The most common and effective model to study streamer dynamics is 
based on the following drift-diffusion equations for electrons and ions, 
coupled with Poisson's equation \cite{Babaeva:1996,Kulikovsky:1997c}:
%%%%%%%%
\def\SeP{\ne\alpha |\vec v_{\rm e}|}
\def\SpP{\ne\alpha |\vec v_{\rm e}|}
\def\SeM{\ne\eta  |\vec v_{\rm e}| +  \ne\np\beta_{\rm ep}}
\def\SpM{\ne\np\beta_{\rm ep} + \nn\np\beta_{\rm np}}
\def\SnM{\nn\np\beta_{\rm np}}
\def\SnP{\ne\eta  |\vec v_{\rm e}|}
\def\SD{\nn\gamma}
%%%%%%%%%
\begin{equation}\label{trasp} 
\left.
\begin{array}{rcl}
\partial_t \ne -\partial_ \vec{x}\cdot\ne\,\ve 
  -\partial_ \vec{x}\cdot(\De\ \partial_ \vec{x}\ne) &=& \SeP -\SeM + \SD,\\ 
\partial_t \np +\partial_ \vec{x}\cdot\np\vp 
  -\partial_ \vec{x}\cdot(\Dp\,\partial_ \vec{x}\np) &=& \SpP -\SpM,\\ 
\partial_t \nn -\partial_ \vec{x}\cdot\nn\vn 
  -\partial_ \vec{x}\cdot(\Dn\,\partial_ \vec{x}\nn)  &=& \SnP -\SnM - \SD, 
\end{array}
\right\}
\end{equation}
\begin{equation}
\varepsilon_0\, \partial_ \vec{x}^2 V = -q_{\rm e}(\np-\nn-\ne), \label{poisson}
\end{equation}
where 
$\vec{x}\in \mathbb{R}^d$,
$n_i$ is the  density of species $i$ (e: electrons, p: positive ions, n: negative ions),
$V$ is the electric potential,
$\vec{v}_i= \mu_i \vec E$ ($\vec E$ being the electric field) is the drift velocity.
$D_i$ and $\mu_i$, are the diffusion coefficient and the absolute value of
mobility of the charged species $i$,
$q_{\rm e}$ is the absolute value of an electron charge, 
and   $\varepsilon_0$  is the permittivity of free space.
$\alpha$ is the impact ionization coefficient, $\eta$
stands for the electron attachment on neutral molecules, $\beta_{\rm ep}$
and $\beta_{\rm np}$
account respectively for the electron-positive ion
and the
negative-positive ion recombination,
and $\gamma$ is the detachment coefficient.

The electric field $\vec E$ and the potential $V$ are related by
\begin{equation}
\vec E = - \partial_ \vec{x} V,
\end{equation}
and thus, the Poisson's equation (\ref{poisson}) becomes:
\begin{equation}\label{poisson2}
\varepsilon_0\, \partial_ \vec{x} \cdot \vec E = q_{\rm e}(\np-\nn-\ne).
\end{equation}

All the coefficients of the model are assumed to
be functions of the local reduced electric field $E/N_{\rm gas}$,
where $E$ is the electric field magnitude and $N_{\rm gas}$ is the
air neutral density. For test studies presented in this
paper, the transport parameters 
for air are taken from \cite{Morrow:1997};
detachment and attachment coefficients,
respectively from \cite{Benilov:2003} and \cite{Kossyi:1992};
and other reaction rates, also from \cite{Morrow:1997}.
Diffusion coefficients for ions are derived from mobilities using 
classical Einstein relations.

In simulations of positive streamer discharges in air at atmospheric pressure without any preionization, 
the photoionization term is crucial to produce seed charges in front of the streamer head and then 
to ensure the streamer propagation \cite{Bourdon:2007}.  
However, in repetitive discharges, \cite{Pancheshnyi:2005b} 
and recently \cite{Wormeester:2010} have shown that even at low frequency, a significant amount
of seed charges from previous discharges may be present in the inter-electrode gap. 
In this work, 
we have neglected the photoionization source term and considered discharge conditions with 
a preionization background to ensure a stable propagation of
the discharge without impacting the main discharge characteristics 
\cite{Pancheshnyi:2001,Pancheshnyi:2005b,Bourdon:2010,Celestin_thesis}.

\section{Construction of the numerical strategy} \label{sec:num_str}
In this section, we introduce a new numerical technique
for multi-scale streamer discharge simulations,
based on a second order decoupled resolution of the electric field and the drift-diffusion
equations for electrons and ions, with self-adaptive decoupling
time steps with error control.
The drift-diffusion equations are then solved 
using a dedicated Strang time operator splitting scheme
for multi-scale phenomena.
On the other hand,
the electric field is computed based on a
parallel computing method, specially conceived for the
configuration under study in 1.5D geometry.
Both resolutions are conducted 
on a dynamic adaptive mesh using
spatial multiresolution transformation with error control of
the spatial adapted representation.

\subsection{Second order adaptive time integration strategy}
\label{splitting_ef}
Let us write the semi-discretized equations 
(\ref{trasp}) and
(\ref{poisson2}) in the following way just for analysis purposes:
\begin{equation}\label{semi_disc} 
\left.
\begin{array}{rcl}
{\rm d}_t \psi &=& \Psi(\psi,\phi), \\
0&=& \Phi(\psi,\phi),
\end{array}
\right\}
\end{equation}
for $t>t_0$, where $\psi:\mathbb{R}\to\mathbb{R}^{N \times m}$
and 
$\phi:\mathbb{R}\to\mathbb{R}^{N \times d}$
stand respectively for the spatial discretization
of $(\ne,\np,\nn)$, i.e. $m=3$,  
and of $\vec E$ over $N$ points.
Supposing that all functions are sufficiently 
differentiable in all their variables
and using the Taylor expansion of the true solution,
one can write after some time $\Delta t$ from initial time $t_0$,
\begin{equation}\label{taylor_ex}
\psi(t_0+\Delta t) = \psi_0 + \Delta t \Psi(\psi_0,\phi_0)
+ \displaystyle{\frac{\Delta t^2}{2} } \left[\partial_\psi \Psi\, \Psi+  \partial_\phi \Psi\, {\rm d}_t \phi \right]_{t=t_0} 
+ \mathcal{O}(\Delta t^3),
\end{equation}
with $\psi_0=\psi(t_0)$, $\phi_0=\phi(t_0)$.

A second order in time resolution of system (\ref{semi_disc})
must then verify (\ref{taylor_ex}) locally for each $\Delta t$.
However, as it was stated before,
solving simultaneously (\ref{trasp}) and (\ref{poisson}) (or (\ref{poisson2})), 
or equivalently (\ref{semi_disc}),
involves important numerical difficulties, 
considering for instance
the different nature of equations (\ref{trasp}) and
(\ref{poisson}) (or (\ref{poisson2})).
Therefore, a decoupled approach is often used in which one aims at solving
the drift-diffusion equations and the electric field independently.
This amounts to
solve 
\begin{equation}\label{semi_disc_t0} 
{\rm d}_t \tilde{\psi} = \Psi(\tilde{\psi},\phi^\star), \quad
 t \in \,]t_0,t_0+\Delta t],
\end{equation}
with fixed $\phi^\star=\phi(t^\star)$, $t^\star \in [t_0,t_0+\Delta t]$ and $\tilde{\psi}(t_0)=\psi_0$.

The most common technique considers $t^\star=t_0$, that is,
to previously compute the electric field at $t_0$ from $\Phi(\psi_0,\phi_0)=0$,
and then solve (\ref{semi_disc_t0}) with $\phi^\star=\phi_0$.
This can be interpreted as a standard first order operator splitting method
that yields an
approximation of order 1, $\tilde{\psi}_1(t)$, of the exact solution, $\psi(t)$,
based on classical numerical analysis results
obtained by
confronting (\ref{taylor_ex}) with
\begin{equation}\label{taylor_lie}
\tilde{\psi}_1(t_0+\Delta t) = \psi_0 + \Delta t \Psi(\psi_0,\phi_0)
+ \displaystyle{\frac{\Delta t^2}{2} } \left[\partial_\psi \Psi\, \Psi \right]_{t=t_0} 
+ \mathcal{O}(\Delta t^3).
\end{equation}
The same result follows for $\tilde{\phi}_1(t_0+\Delta t)$ computed out of 
$\Phi(\tilde{\psi}_1(t_0+\Delta t),\tilde{\phi}_1(t_0+\Delta t))=0$ or
equivalently, out of its explicit representation 
$\tilde{\phi}_1(t_0+\Delta t)= \Upsilon (\tilde{\psi}_1(t_0+\Delta t))$,
assuming a Lipschitz condition:
\begin{equation}
\left\| \Upsilon (\psi)-\Upsilon (\psi^\star)\right\|
\leq L \left\|\psi-\psi^\star\right\|.
\end{equation}
Considering now any $t^\star \in [t_0,t_0+\Delta t]$ into (\ref{semi_disc_t0}), the only 
second order solution, $(\tilde{\psi}_2(t),\tilde{\phi}_2(t))$, will be given by 
resolution of (\ref{semi_disc_t0})
with $\phi^\star=\phi_{\frac{1}{2}}$ for $t^\star=t_0+\Delta t/2$, for which
\begin{equation}\label{taylor_strang1}
\tilde{\psi}_2(t_0+\Delta t) = \psi_0 + \Delta t \Psi(\psi_0,\phi_{\frac{1}{2}})
+ \displaystyle{\frac{\Delta t^2}{2} } \left[\partial_\psi \Psi\, \Psi \right]_{\psi \to \psi_0,\phi_{\frac{1}{2}}} 
+ \mathcal{O}(\Delta t^3),
\end{equation}
where
\begin{eqnarray}
\Psi(\psi_0,\phi_{\frac{1}{2}}) &=&
\Psi\left(\psi_0,\phi\left(t_0+\frac{\Delta t}{2}\right)\right)\nonumber\\
&=&
\Psi\left(\psi_0,\phi_0 + \displaystyle{\frac{\Delta t}{2} } \left.{\rm d}_t \phi \right|_{t=t_0}+
\mathcal{O}(\Delta t^2)\right)\nonumber\\
&=&
\Psi(\psi_0,\phi_0)+
\displaystyle{\frac{\Delta t}{2} }
\left[\partial_\phi \Psi\, {\rm d}_t \phi \right]_{t=t_0}
+\mathcal{O}(\Delta t^2),
\end{eqnarray}
and hence,
\begin{equation}\label{taylor_strang2}
\tilde{\psi}_2(t_0+\Delta t) = \psi_0 + \Delta t \Psi(\psi_0,\phi_0)
+ \displaystyle{\frac{\Delta t^2}{2} } \left[\partial_\psi \Psi\, \Psi+  \partial_\phi \Psi\, {\rm d}_t \phi \right]_{t=t_0} 
+ \mathcal{O}(\Delta t^3);
\end{equation}
and 
\begin{equation}
\tilde{\phi}_2(t_0+\Delta t)= \Upsilon (\tilde{\psi}_2(t_0+\Delta t)).
\end{equation}

Nevertheless, this second order approximation, $\tilde{\psi}_2(t)$,
is based on the previous knowledge of $\phi_{\frac{1}{2}}=\phi(t_0+\Delta t/2)$,
and thus, of $\psi(t_0+\Delta t/2)$.
In order to overcome this difficulty, one can solve
(\ref{semi_disc_t0}) with $\phi^\star=\tilde{\phi}_1(t_0+\Delta t/2)= \Upsilon (\tilde{\psi}_1(t_0+\Delta t/2))$,
that is, computing first $\tilde{\psi}_1(t_0+\Delta t/2)$ with the first order
method.
In particular, this does not change the previous order estimates 
as it follows from
\begin{eqnarray}
\psi(t_0+\Delta t)-\tilde{\psi}_2(t_0+\Delta t)&=&
\displaystyle{\frac{\Delta t^2}{2} } \left[\partial_\phi \Psi\, {\rm d}_t (\phi-\tilde{\phi}_1) \right]_{t=t_0}
+ \mathcal{O}(\Delta t^3)\nonumber\\
&=&
\displaystyle{\frac{\Delta t^2}{2} } \left[\partial_\phi \Psi\, \partial_\psi \Upsilon\, {\rm d}_t (\psi-\tilde{\psi}_1) \right]_{t=t_0}
+ \mathcal{O}(\Delta t^3)\nonumber\\[1ex]
&=& \mathcal{O}(\Delta t^3).
\end{eqnarray}

Taking into account both methods,
\begin{equation}
\left(
\begin{array}{c}
 \tilde{\psi}_1(t_0+\Delta t)\\
\tilde{\phi}_1(t_0+\Delta t)
\end{array}
\right)
=
\mathcal{T}_1^{\Delta t}
\left(
\begin{array}{c}
 \psi_0\\
\phi_0
\end{array}
\right),
\quad
\left(
\begin{array}{c}
 \tilde{\psi}_2(t_0+\Delta t)\\
\tilde{\phi}_2(t_0+\Delta t)
\end{array}
\right)
=
\mathcal{T}_2^{\Delta t}
\left(
\begin{array}{c}
 \psi_0\\
\phi_0
\end{array}
\right),
\end{equation}
we perform computations with a
second order scheme
$\mathcal{T}_2^{\Delta t}$,
which uses 
an embedded and lower order scheme
$\mathcal{T}_1^{\Delta t/2}$,
as it was previously detailed.
An adaptive time step strategy
is then implemented in order to 
control the accuracy of computations
by tuning the duration of the
decoupled resolution.
It is based on a local
error estimate, dynamically computed
at the end of each decoupling time step $\Delta t$,  
given by 
\begin{equation}\label{est_split_ef}
 \big\|\mathcal{T}_2^{\Delta t}(\psi_0,\phi_0)^t
-\mathcal{T}_1^{\Delta t}(\psi_0,\phi_0)^t\big\| 
\approx \mathcal{O}({\Delta t}^2).
\end{equation}

Therefore,
for a given accuracy tolerance $\eta_{\mathcal{T}}$,
\begin{equation}\label{est_split_er_tol}
 \big\|\mathcal{T}_2^{\Delta t}(\psi_0,\phi_0)^t
-\mathcal{T}_1^{\Delta t}(\psi_0,\phi_0)^t\big\| 
< \eta_{\mathcal{T}}
\end{equation}
must be verified
in order to accept the current computation with $\Delta t$,
while the new time step is calculated by
\begin{equation}\label{delta_split_ef}
\Delta t^{\rm new} = \Delta t  
\sqrt{\frac{\eta_{\mathcal{T}}}
{\big\|\mathcal{T}_2^{\Delta t}(\psi_0,\phi_0)^t
-\mathcal{T}_1^{\Delta t}(\psi_0,\phi_0)^t\big\|}}.
\end{equation}

Several dedicated solvers can be then implemented for each subproblem
(\ref{trasp}) and (\ref{poisson})
while the theoretical error estimates 
of the decoupling schemes
analyzed in this section
remain valid.
In this way,
the independent choice of appropriate numerical schemes
allows
to strongly reduce the computational complexity of the 
global numerical strategy,
and an error control procedure such as the one 
proposed in this work
allows to effectively calibrate this decoupling within a
prescribed accuracy tolerance.

\subsection{Resolution of the drift-diffusion equations}
We consider now the numerical resolution of the drift-diffusion equations (\ref{trasp}),
that can be written in the general form of a convection-reaction-diffusion system of equations: 
\begin{equation}\label{sys_reac_diff_gen}
\left.
\begin{array}{rcll}
\displaystyle{
\partial _t \vec{u} - \partial_{\vec{x}} \left( 
\vec{F} \left( \vec{u} \right) +
\vec{D}(\vec{u})
 \partial_{\vec{x}} \vec{u} \right)
}
&=&
{\displaystyle
\vec{f} 
\left(\vec{u}\right),
}
& t>t_0,
\\
\displaystyle{
\vec{u}(t_0,\vec{x})
}
&=&
{\displaystyle
\vec{u}_0(\vec{x}),
}
& t=t_0,
\end{array}
\right\}
\end{equation}
where $\vec{F}$, $\vec{f}:\mathbb{R}^{m}\to\mathbb{R}^{m}$ and 
$\vec{u}:\mathbb{R}\times \mathbb{R} ^d \to\mathbb{R}^{m}$,
with a tensor of order $d\times d\times m$ as
diffusion matrix $ \vec{D}(\vec{u})$.
In particular, $\vec{u} = (\ne,\np,\nn)^t$ with $m=3$ in this study.

The system (\ref{sys_reac_diff_gen}) 
corresponds to problem (\ref{semi_disc_t0}) for a fixed electric field,
and
it is solved during each decoupling time step $\Delta t$
into $\mathcal{T}_2$ (or $\mathcal{T}_1$) scheme,
using a
Strang time operator scheme with dedicated high order 
time integrators on a dynamic adaptive mesh, based on a strategy introduced
in \cite{article_mr}.
This resolution is briefly detailed in following sections.

\subsubsection{Time operator splitting}
\label{splitting_strang}
An operator splitting procedure allows to consider dedicated
solvers for 
the reaction part 
which is decoupled
from other physical phenomena
like convection, diffusion or both, for which there also exist dedicated
numerical methods. 
These dedicated methods chosen for
each subsystem are then responsible for dealing with the fastest scales associated
with each one of them, in a separate manner,
while the reconstruction of the global solution by
the splitting scheme should guarantee
an accurate description with error control of the global
physical coupling, without being related to the stability constraints of the numerical resolution of each subsystem.

Considering problem (\ref{sys_reac_diff_gen}) and in order to remain consistent
with the second order $\mathcal{T}_2$ scheme, a 
second order Strang scheme is implemented \cite{Strang68}
\begin{equation}\label{strang}
\mathcal{S}^{\Delta t_{\rm s}}(\vec{u}_0) 
= 
\mathcal{R}^{\Delta t_{\rm s}/2}
\mathcal{D}^{\Delta t_{\rm s}/2}
\mathcal{C}^{\Delta t_{\rm s}}
\mathcal{D}^{\Delta t_{\rm s}/2}
\mathcal{R}^{\Delta t_{\rm s}/2}
(\vec{u}_0),
\end{equation}
where operators $\mathcal{R}$, $\mathcal{D}$, $\mathcal{C}$ 
indicate respectively the independent resolution of
the reaction, diffusion and convection problems 
with a splitting time step,
$\Delta t_{\rm s}$, 
taken inside the overall decoupling time step,
$\Delta t_{\rm s}\leq \Delta t$.
Usually,
for propagating reaction waves
where for instance,
the speed of propagation is much slower than some of the chemical
scales,
the fastest scales are not
directly related to the global physics of the phenomenon, 
and thus,
larger splitting time steps might be considered
\cite{article_mr,article_avc}.
Nevertheless,
order reductions may then appear due to short-life transients 
associated with fast variables and
in these cases,
it has been proven in \cite{Descombes04} 
that better performances are expected while ending
the splitting scheme by operator $\mathcal{R}$
or in a more general case,
the part involving the fastest time scales of the phenomenon.

The resolution of (\ref{sys_reac_diff_gen}) should be precise enough
to guarantee theoretical estimates given in Section~\ref{splitting_ef}.
Therefore, an adaptive splitting time step strategy, based on a local
error estimate at the end of each splitting time step $\Delta t_{\rm s}$, 
is also implemented in order to 
control the accuracy of computations \cite{article_fvca}.
In this context,
a second, embedded and lower order Strang splitting method 
$\widetilde{\mathcal{S}}^{\Delta t_{\rm s}}$ was developed 
by \cite{article_pas},
which allows to dynamically calculate a local error estimate
that should verify
\begin{equation}\label{est_split}
 \big\|\mathcal{S}^{\Delta t_{\rm s}}(\vec{u}_0)
-\widetilde{\mathcal{S}}^{\Delta t_{\rm s}}(\vec{u}_0)
\big\| \approx \mathcal{O}({\Delta t_{\rm s}}^2) 
< \eta_{\rm split},
\end{equation}
in order to accept the current computation with $\Delta t_{\rm s}$,
and thus, the new splitting time step is given by
\begin{equation}\label{delta_split}
\Delta t_{\rm s}^{\rm new} = 
\min \left( 
\Delta t_{\rm s}  \sqrt{\frac{\eta_{\rm split}}{\big\|\mathcal{S}^{\Delta t_{\rm s}}(\vec{u}_0)
-\widetilde{\mathcal{S}}^{\Delta t_{\rm s}}(\vec{u}_0)\big\|}},\ 
t_0+\Delta t - \hat{t}
\right),
\end{equation}
with $\eta_{\rm split}\leq \eta_{\mathcal{T}}$ and $\hat{t} =\sum_i \Delta t_{{\rm s}_i}$
while $\hat{t} \in\, ]t_0,t_0+\Delta t]$.

The choice of suitable time integration methods to 
numerically approximate  
$\mathcal{R}$, $\mathcal{D}$ and $\mathcal{C}$
during each $\Delta t_{\rm s}$
is mandatory not only to guarantee the theoretical framework
of the numerical analysis
but also to take advantage of the particular features
of each independent subproblem.
A new operator splitting for reaction-diffusion systems was 
recently introduced \cite{article_mr,article_avc},
which considers
a high fifth order, \textit{A}-stable,
\textit{L}-stable method like Radau5 \cite{Hairer96},
based on implicit Runge-Kutta schemes for stiff ODEs, 
that solves
with a local cell by cell approach
the reaction term: a system of stiff ODEs without spatial coupling
in a splitting context.
For the diffusion problem, another high fourth order method like 
ROCK4 \cite{Abdulle02} is considered, 
which is based on explicit stabilized Runge-Kutta schemes
that feature extended
stability domains along the negative real axis.
The ROCK4 solver is then very appropriate for diffusion
problems because of the usual 
predominance of negative real eigenvalues.
Both methods incorporate adaptive time integration tools,
similar to (\ref{delta_split_ef}) and (\ref{delta_split}),
in order to control the accuracy of the integrations for given 
accuracy tolerances
$\eta_{\rm Radau5}$ and
$\eta_{\rm ROCK4}$, chosen such that
$\eta_{\rm Radau5}<\eta_{\rm split}$ and $\eta_{\rm ROCK4}<\eta_{\rm split}$.
In particular,
in the case of multi-scale propagating waves,
it can be proven that 
the local treatment plus the adaptive time stepping 
of the reaction solver
allow to discriminate the cells of high reactive activity 
only present in the neighborhood of the localized wavefront, 
saving as a consequence a large quantity of integration time \cite{article_avc}.

An explicit high order in time and in space one step monotonicity 
preserving scheme OSMP \cite{daru04} is used
as convective scheme. It combines monotonicity preserving
constraints for non-monotone data to avoid extrema clipping,
with TVD features to prevent spurious oscillations around 
discontinuities or sharp spatial gradients.
Classical CFL stability restrictions are though imposed
inside each splitting time step $\Delta t_{\rm s}$ for operator 
$\mathcal{C}^{\Delta t_{\rm s}}$.
The overall combination of an explicit treatment of
the spatial phenomena as convection and diffusion, with a local
implicit integration of stiff reaction implies important
savings in computing time and memory resources \cite{article_mr},
as well as an important reduction of computational
complexity with respect to a fully implicit coupled resolution 
of problem (\ref{sys_reac_diff_gen}).
On the other hand an explicit coupled treatment of (\ref{sys_reac_diff_gen})
will have a very
limited efficiency for stiff problems unless more sophisticated 
strategies as the asynchronous
local time-stepping techniques 
\cite{Karimabadi05,Omelchenko06,Unfer07}
are considered even though these schemes 
do not provide a precise measurement of the accuracy of
the integration.

Finally, 
the numerical errors of the splitting scheme are effectively handled 
by an error control procedure which furthermore allows to determine
the coupling time scales of the global phenomenon that can be
several orders of magnitude slower than the fastest time scales
of each subproblem treated by each dedicated solver.
In this way a decoupling of the time scale spectrum of the problem
is achieved that leads to more efficient performances within a prescribed
accuracy tolerance whenever this decomposition of scales is possible.

\subsubsection{Mesh refinement technique}
Regarding problem (\ref{sys_reac_diff_gen}),
we are concerned with
the propagation of reacting wavefronts, hence important 
reactive activity as well as steep spatial gradients
are localized phenomena.
This implies that if we consider the resolution of
the reaction problem, a con\-si\-de\-ra\-ble
amount of computing time is spent on nodes
that are practically at (partial) equilibrium.
Moreover, there is no need to represent these quasi-stationary
regions with the same spatial discretization needed to describe
the reaction front, 
so that the convection and the diffusion problems
might also be solved over a smaller number of nodes. 
An adapted mesh obtained by a mul\-ti\-re\-so\-lu\-tion
process  which discriminates the various space scales 
of the phenomenon,
turns out to be a very convenient solution to overcome
these difficulties
\cite{article_mr,article_esaim}.
Furthermore, 
in plasma applications, the resolution of Poisson's equation takes
usually $\sim$$80\%$ of the computing time.
Thus, important savings are achieved with a mesh adaptive technique,
as a consequence of the strong reduction of cells.

In practice,
if one considers a set of nested spatial grids
from the coarsest to the finest one,
a multiresolution transformation allows to 
represent a discretized function as values on the
coarsest grid plus a series of local estimates  
at all other levels of such nested grids \cite{Cohen03}.
These estimates 
correspond to the wavelet coefficients of a wavelet decomposition
obtained by inter-level 
transformations,
and retain 
the information on local regularity when
going from a coarse to a finer grid.
Hence, the main idea is to use the decay of the
wavelet coefficients to obtain information on the local 
regularity of the solution: 
lower wavelet coefficients are associated
with locally regular spatial configurations and vice-versa.
The basis of this strategy is presented in the following.
For further details 
on adaptive
multiresolution techniques, we refer to the books of 
\cite{Cohen:2000} and 
\cite{Muller:2003}.

\subsubsection{Basis of a multiresolution representation}
\label{sec_mr}
To simplify the presentation
let us consider 
nested finite volume discretizations of 
(\ref{sys_reac_diff_gen})
with only one component, $m=1$.
For $l=0,1,\cdots,L$ 
from the coarsest to the finest grid,
we have then regular disjoint partitions (cells) $(\Omega_{\gamma})_{\gamma \in S_l}$
of an open subset $\Omega \subset \mathbb{R}^{d}$,
such that each $\Omega _{\gamma}$, $\gamma \in S_l$, 
is the union of a finite number of cells $\Omega_{\mu}$, $\mu \in S_{l+1}$,
and thus, $S_l$ and $S_{l+1}$ are consecutive embedded grids.
We denote
$\mathbf{U}_l:=(u_{\gamma})_{\gamma \in S_l}$ as the 
representation of $u$
on the grid $S_l$ where $u_{\gamma}$
represents the cell-average of 
$u :\mathbb{R}\times \mathbb{R}^d \to\mathbb{R}$
in $\Omega _{\gamma}$,
\begin{equation}
  u_{\gamma} := |\Omega_{\gamma}|^{-1} \int_{\Omega_{\gamma}}  u(t,\mathbf{x})\, {\rm d}\mathbf{x}.
 \end{equation}

The data at different levels of discretization are related
by two inter-level transformations which are defined as follows:
\begin{enumerate}
 \item 
The \emph{projection} operator $P^l_{l-1}$, 
which maps $\mathbf{U}_l$ to $\mathbf{U}_{l-1}$.
It is obtained
through exact a\-ve\-ra\-ges computed at the finer level by
\begin{equation}\label{average}
  u_{\gamma} = |\Omega_{\gamma}|^{-1} 
\sum_{|\mu|=|\gamma| +1,\Omega_{\mu} 
\subset \Omega_{\gamma}} |\Omega_{\mu}|u_{\mu},
 \end{equation}
where $|\gamma|:=l$ if $\gamma \in S_l$.
As far as grids are nested, this
projection operator is \emph{exact} and \emph{unique}
\cite{Cohen:2000}.
\item
The \emph{prediction} operator $P^{l-1}_l$,
which maps $\mathbf{U}_{l-1}$ to 
an approximation $\hat{\mathbf{U}}_l$ of $\mathbf{U}_{l}$.
There is an
infinite number of choices to define $P^{l-1}_l$, but we impose at least two basic
constraints \cite{Cohen03}:
\begin{enumerate}
\item The prediction is local, i.e., $\hat{u}_{\mu}$ for a given $\Omega_{\mu}$ 
depends on a set of values $u_{\gamma}$
 in a finite stencil $R_{\mu}$ surrounding $\Omega_{\mu}$, where $|\mu|=|\gamma| +1$.
\item The prediction is consistent with the projection in the sense that
 \begin{equation}\label{consis}
   |\Omega_{\gamma}| u_{\gamma} = \sum_{|\mu|=|\gamma| +1,\Omega_{\mu} \subset \Omega_{\gamma}} |\Omega_{\mu}|\hat{u}_{\mu};
 \end{equation}
    i.e., $P_{l-1}^l \circ P_l^{l-1} = Id$. 
 \end{enumerate}
\end{enumerate}
With these operators, we define for each cell $\Omega _{\mu}$
the prediction error or \textit{detail} 
as the difference between the exact and predicted values:
\begin{equation}\label{detail}
d_{\mu} := u_{\mu} - \hat{u}_{\mu},
\end{equation}
or in terms of inter-level operations:
\begin{equation}\label{detail2}
d_{\mu} = u_{\mu} - P_{|\mu|}^{|\mu|-1} \circ P_{|\mu|-1}^{|\mu|} u_{\mu}.
\end{equation}
We can then construct a \textit{detail vector} $\mathbf{D}_l$
as shown in \cite{Cohen03} in order to 
get a one-to-one correspondence 
from expressions
(\ref{detail}) and (\ref{consis}): 
\begin{equation}
\mathbf{U}_l\longleftrightarrow (\mathbf{U}_{l-1},\mathbf{D}_l).
\end{equation}
Hence, 
by iteration of this decomposition, we finally obtain a multi-scale representation of $\mathbf{U}_L$
in terms of $\mathbf{M}_L = (\mathbf{U}_0,\mathbf{D}_1,\mathbf{D}_2,\cdots,\mathbf{D}_L)$:
\begin{equation}
\mathcal{M}:\mathbf{U}_L\longmapsto \mathbf{M}_L,
\end{equation}
where the \textit{details} computed with (\ref{detail2}) stand for the wavelet coefficients in a wavelet basis.

One of the main interests of carrying on such a
wavelet decomposition is that this new representation
defines a whole set of regularity
estimators all over the spatial domain and
thus, a data compression might be achieved by deleting
cells whose \textit{detail} verifies
\begin{equation}\label{lambda}
|d_{\mu}| < \varepsilon_l, 
\quad l=|\mu|,
\quad \varepsilon_l = 2^{\frac{d}{2}(l-L)}\eta_{\rm MR},
\end{equation}
where $\eta_{\rm MR}$ is a threshold value defined 
for the finest level $L$ \cite{Harten95}.

An important theoretical result is that
if we denote by $\mathbf{V}^n_L:=(v_{\lambda}^{n})_{\lambda \in S_L}$, 
the solution fully computed
on the finest grid,
and
denote by $\mathbf{U}^n_L$, the solution 
reconstructed on the finest grid
that used
adaptive multiresolution 
(keeping in mind that the time integration was really performed
on a compressed representation of $\mathbf{U}^n$);
then, 
for a fixed time $T=n\Delta t$,
it can be shown \cite{Harten95,Cohen03} that 
the approximation error made
by using this space adaptive representation is proportional to
the threshold value $\eta_{\rm MR}$:
\begin{equation}\label{err_mr}
 \|\mathbf{U}^n_L- \mathbf{V}^n_L \|_{L^2} \propto 
n\eta_{\rm MR}.
\end{equation}

\subsection{Computation of the electric field}
In this part, we are concerned with the resolution of the electric field 
according to the $\mathcal{T}_2$ (or $\mathcal{T}_1$)
scheme at some fixed time for
a given distribution of charges $(\ne,\np,\nn)$,
considering a 1.5D model.
This computation is also performed on the adapted mesh obtained
by the previous multiresolution analysis.

\subsubsection{Discretization of the computational domain}
\begin{figure}[!htb]
 \begin{center}
   \includegraphics[width=0.5\hsize]{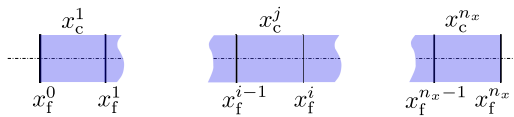}
   \caption{Definition of the grid: the cell centers are located at $\xc^j$, 
whereas cell
faces are located at $\xf^{i}$. The domain is bounded by faces $\xf^0$ (cathode) and $\xf^{n_x}$
(tip of the anode).
 \label{fgrid}}
 \end{center}
\end{figure}

According to Figure~\ref{fdomain},
the computational domain is limited by a planar cathode at $x=0$
and the tip of a hyperbolic anode at $x=L_x$.
The anode is not included in the domain. We consider streamers of fixed radius
$\Rs$  along the axis of symmetry.
The computational domain is divided into $n_x$ cells
of different size corresponding to the multiresolution adapted mesh,
 with faces $\xf^i$, where
$i\in[0,n_x]$
and  cell centers $\xc^j$, where $j\in [1,n_x]$. The face $\xf^0$ corresponds
to the position of the cathode and $\xf^{n_x}$ corresponds to the position
of the tip of the anode. Therefore for each cell $\xc^i$, there is its left
face
$\xf^{i-1}$, and its right face $\xf^{i}$. For each cell $\xc^j$ we define
a width
$w_j=\xf^j-\xf^{j-1}$ (see Figure~\ref{fgrid}).

\subsubsection{Resolution of the electric field in a 1.5D model}
To determine the electric field during the propagation of the streamer, the space 
charge of the streamer is considered as a set of finite cylinders of width $w_j$, bounded by
cell faces $x_{\rm f}^{j-1}$ and $x_{\rm f}^{j}$. 
As the computational domain is bounded by conducting electrodes of fixed potential,
each volume charge $\rho_j$ creates an infinite series of image charges \cite{Davies:1964,Davies:1971}. 
Then the principle of superposition is used to sum individual contributions from
all the cylindrical
space charges in the domain, their image charges, and the Laplacian electric field
(computed based on classical results \cite{Eyring28}). An advantage of this approach
dwells in the fact that the electric field contributions from individual cylinders
can be expressed analytically in a simple form and the determination of the 
electric field in each point of the domain can be performed in parallel.

\begin{figure}[!htb]
 \begin{center}
   \includegraphics[width=0.85\hsize]{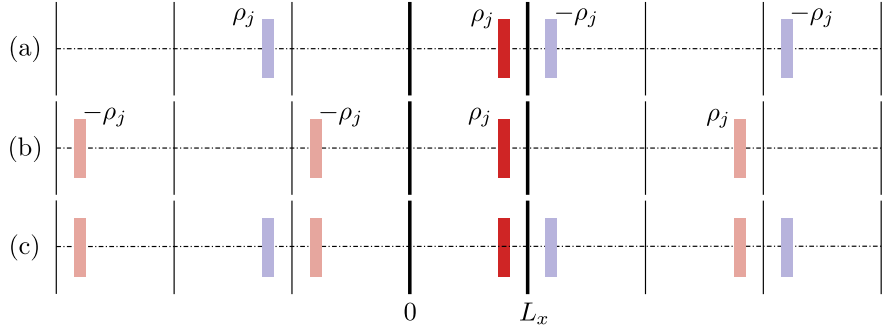}
   \caption{Image charges up to the third order: (a) charge $\rho_j$ is 
            first mirrored behind the anode ($x=L_x$), (b) charge $\rho_j$ is
            first mirrored behind the cathode ($x=0$), (c) charge $\rho_j$ 
            and its images.
            \label{images.eps}
}
 \end{center}
\end{figure}

In the configurations we have studied, the cathode is grounded whereas an electric voltage
is applied on the anode. These boundary conditions are taken into account
by the Laplacian electric field and by including a series of image charges of the charges
in the gap. It is important to note that the computation of 
the Laplacian electric field takes into account the real
geometry of electrodes as shown in Figure~\ref{fdomain}.
However, in this work,
to simplify the computation of image charges
we have assumed that both electrodes are planar.
For a volume charge $\rho_j$  centered  at $\xc^j$,  
there exist image charges of the first order with charge $-\rho_j$ at 
$x=2L_x-\xc^j$ mirrored through the anode, see Figure~\ref{images.eps}a,
and at $x=-\xc^j$ mirrored through the cathode, see 
Figure~\ref{images.eps}b. And for each of these image charges
there exist higher order  image charges of opposite signs and so forth.
All the image charges of $\rho_j$ up to order three are depicted
in Figure~\ref{images.eps}c.

\begin{figure}[!htb]
 \begin{center}
   \includegraphics[width=0.5\hsize]{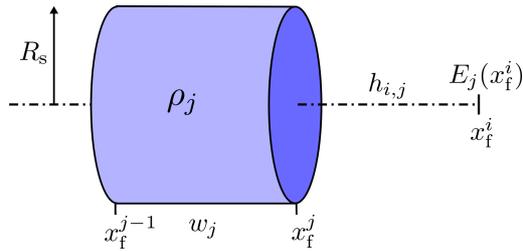}
   \caption{Charged cylinder considered to compute the electric field in the 1.5D model.
 \label{valec.eps}}
 \end{center}
\end{figure}

Integrating the generalized Coulomb's law \cite{Jackson:1999} and using the principle
of superposition, we find that the cylinder charges of cells $j\in[1,n_x]$ 
of width $w_j$, radius $\Rs$, charged with densities $\rho_j$
(see Figure~\ref{valec.eps}),
and the Laplacian electric field $E_L(\xf^i)$ at $\xf^i$ \cite{Eyring28},
create the electric field $E$ at position $\xf^i$ as follows:
\begin{equation}\label{EFi}
E(\xf^i)= E_L(\xf^i) + \sum_{j=1}^{n_x} s{\rho_{j}w_j\over2\epsi}\left(1 -
 {w_j+2h_{i,j}\over\sqrt{h_{i,j}^2 +\Rs^2} +
\sqrt{\left(h_{i,j}+w_j\right)^2 + \Rs^2}} \right),
\end{equation}
where
$$
h_{i,j}=\cases{ \xf^i -\xf^j & for $ i \geq j $  \cr
                        \xf^{j-1} - \xf^i & for $  i < j $}
\quad
\hbox{and}
\quad
s=\cases{ +1 & for $ i \geq j $  \cr
          -1 & for $ i < j $        }
$$

The positive sign of $s$ accounts for the electric field calculated  on the right
from
the  position  of the charged cylinder and vice-versa. The same formula applies for
the image charges,
but an appropriate sign of the charge has to be carefully taken into account according to 
Figure~\ref{images.eps}.
In particular, in a shared memory computing environment, 
a straightforward parallelization is accomplished for 
equation (\ref{EFi}),
in which each core solves successively the electric field on 
one single position $\xf^i$, and
where neither synchronization stages nor data exchange are needed among nodes.

Note that for $\Rs\rightarrow\infty$ (infinite plane charges), equation
(\ref{EFi})
yields the exact electric field for a planar front:
$$E_{\rm inf}= \sum_{j=1}^{n_x} s{\rho_{j}w_j\over2\epsi}.$$
For finite radius $\Rs$ the solution (\ref{EFi}) is valid only on the axis of
the discharge,
but when applied to a discharge of a small radius, the electric field will
vary only negligibly
over the cross section of the discharge.  This approach is expected to be
more accurate
for any finite radius than any discretization of Poisson's equation
\cite{Davies:1964}.

%%%%%%%%%%%%%%%%%%%%%%%%%%%%%%%%%%%%%%%%%%%%%%%%%%%%%%%%%%%%%%%%%%%%%%%%%%%%%%%%%%
\section{Numerical results} \label{sec:num_res}
In this  section, we present some numerical illustrations of the proposed numerical strategy for
the simulations of positive streamers using a 1.5D model in a point-to-plane
geometry. 
First, we will consider a discharge propagation with constant applied voltage
for which different features of the numerical strategy are discussed, e.g., error estimates, 
data compression values and computing time, in order
to properly choose the simulation parameters.
Then, the potential of the method is fully exploited for a more complex configuration
of repetitive discharges generated by high frequency pulsed applied voltages,
followed by a long time scale relaxation,
for which a complete
physical description of the discharge and the post-discharge phases
is achieved.

\subsection{Propagation of a positive streamer with constant applied voltage}\label{sec-exact}
We consider a point-to-plane geometry with a $1\,$cm gap
between the tip of the electrode and the plane, and a constant applied
voltage of $13\,$kV at $x=L_x$. For the following simulations, the discharge is initiated by
placing a neutral plasma cloud with a Gaussian distribution close to
the tip of the anode. The initial distributions of electrons and
ions are then given by
\begin{equation} \label{ini_c}
\left.n_{{\rm e,p}}(x)\right|_{t=0} = n_{\rm max}\exp\left(- (x-c)^2/w^2\right) + n_0, 
\quad  \left.n_{\rm n}(x)\right|_{t=0} = 0,
\end{equation}
where $w=0.027\,$cm, $c=1\,$cm, $n_{\rm max}=10^{14}\,$cm$^{-3}$, and with a 
preionization of $n_0=10^{8}\,$cm$^{-3}$. 
There are no negative ions as initial condition.
The streamer radius is set to $\Rs=0.05\,$cm to have a typical  
electric field magnitude in the streamer head of $120\,$kV/cm
 \cite{Kuli_PRE:1998}.
Homogeneous Neumann boundary conditions were considered for the
drift-diffusion equations.
\begin{figure}[!htb]
 \begin{center}
   \includegraphics[width=0.49\hsize]{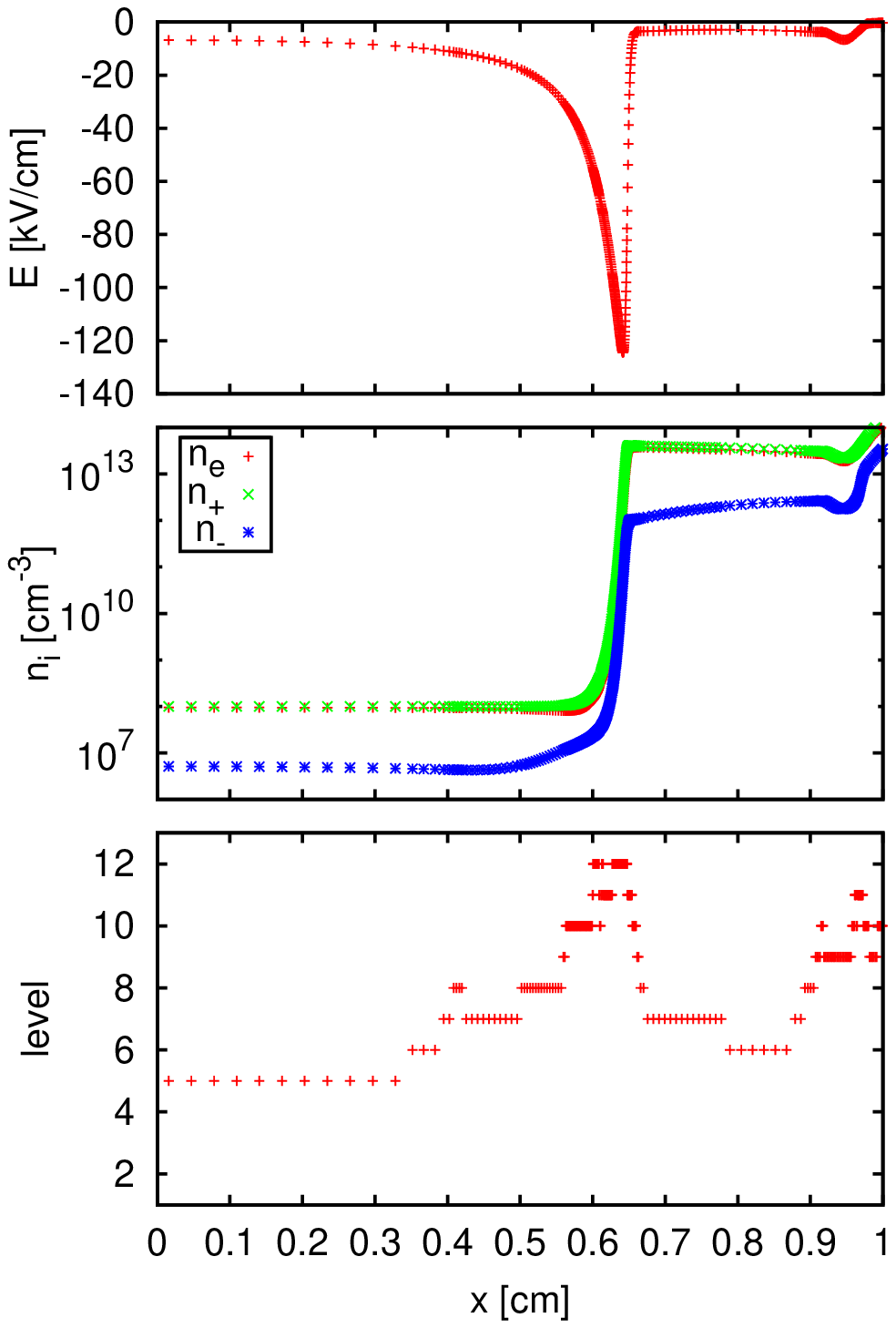} \hfill
   \includegraphics[width=0.49\hsize]{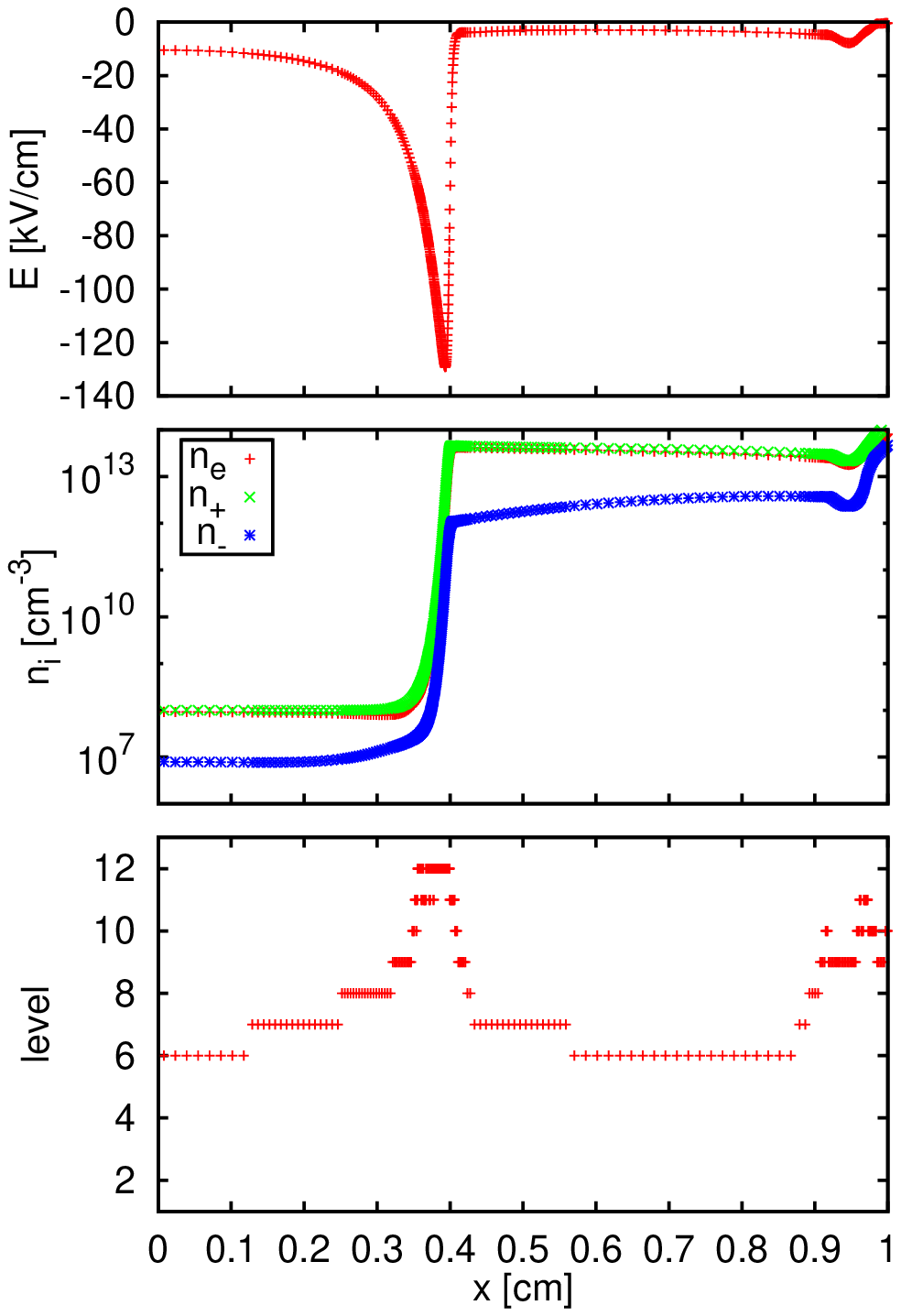} 
   \caption{Positive streamer propagation at $t=6\,$ns (left) and  $t=10\,$ns (right).
Top: electric field; middle: charged species density;
and bottom: grid levels. Finest grid: 4096, $\eta_{\mathcal{T}}=\eta_{\rm split}=\eta_{\rm MR}=10^{-4}$.
 \label{propag}}
 \end{center}
\end{figure}

Two instances of the discharge propagation
are shown in Figure~\ref{propag}, for $12$ nested grids 
equivalent to $4096$ cells on the finest
grid, $L=12$, and for accuracy tolerances of
$\eta_{\mathcal{T}}=\eta_{\rm split}=\eta_{\rm MR}=10^{-4}$;
the spatial refinement takes place only where it is required.
Fine tolerances were chosen in all cases for the 
solvers, $\eta_{\rm Radau5}=\eta_{\rm ROCK4}=10^{-7}$, to guarantee
accurate integrations.
For all the simulation cases, 
the \textit{detail} in each cell is taken as the maximum of the \textit{details} computed according
to (\ref{detail2}) for each variable, where the prediction operator is a 
polynomial interpolation of degree $2$,
performed on normalized $\log \vec{u}$ of the density variables in order to properly
discriminate the streamer heads from the highly
ionized plasma channel;
this logarithmic scale guarantees a correct spatial representation of the phenomenon as seen in
Figure~\ref{propag} for the density profiles.

In order to perform an analysis of the numerical results,
we define the {\it reference} solution
as a fine resolution
with the $\mathcal{T}_2$ scheme that considers a fixed decoupling 
time step, $\Delta t=10^{-14}\,$s and a uniform grid of $4096$ cells.
For this {\it reference} solution,
the memory requirements are acceptable and the simulation is still feasible, but it requires about
14 days of real simulation time on an AMD Opteron 6136 Processor cluster, while running the electric field
computation in parallel on 16 CPU cores. 
In this case, the computation of the electric field,
based on a direct integration of individual contributions
of the charged cylinders,
represents 80\% of total CPU time per time step (about $3.2\,$s). 
\begin{figure}[!htb]
 \begin{center}
   \includegraphics[width=0.49\hsize]{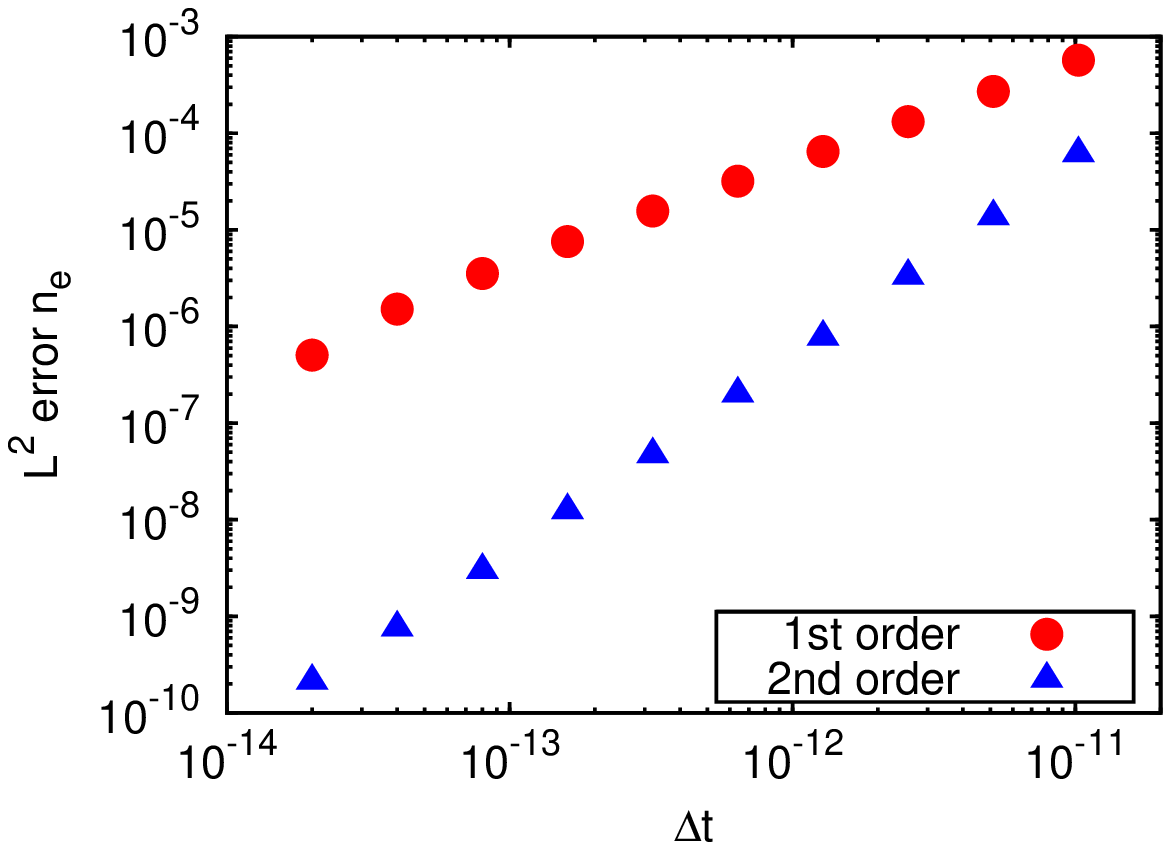}\hfill
  \includegraphics[width=0.49\hsize]{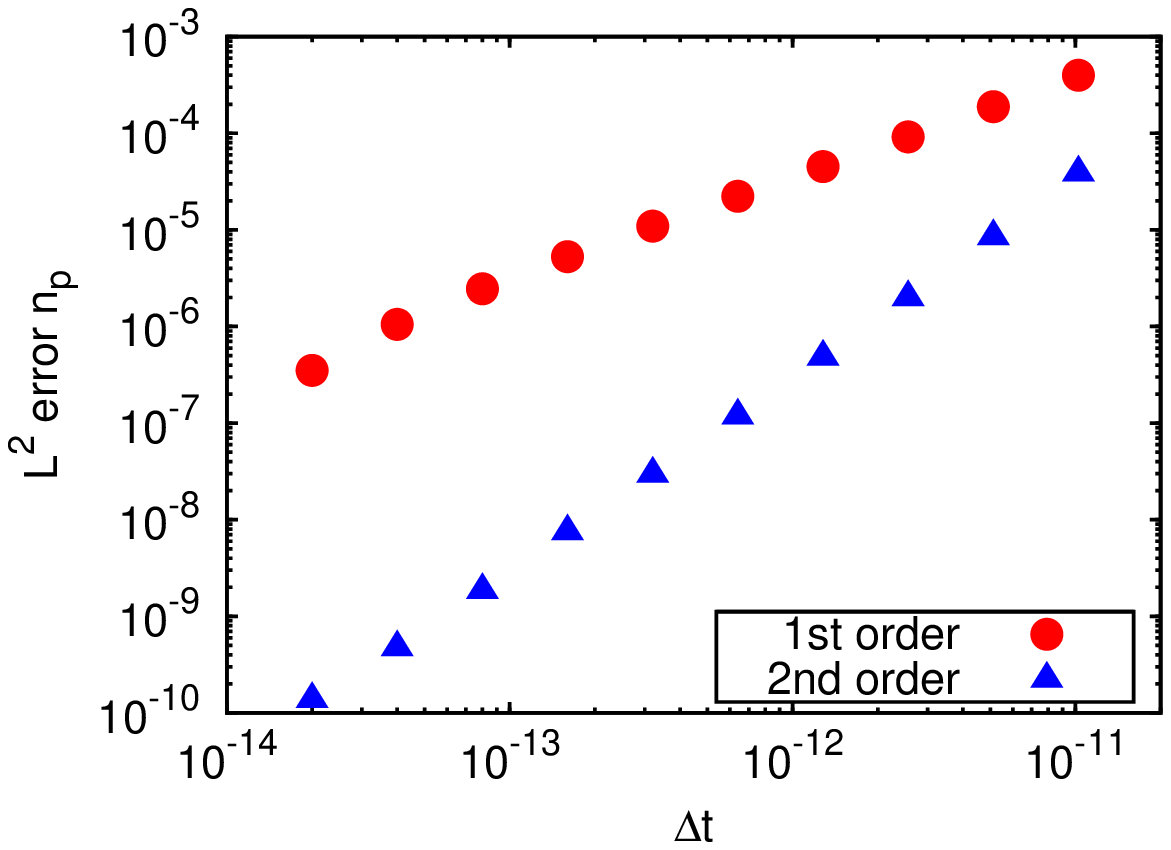}
 \includegraphics[width=0.49\hsize]{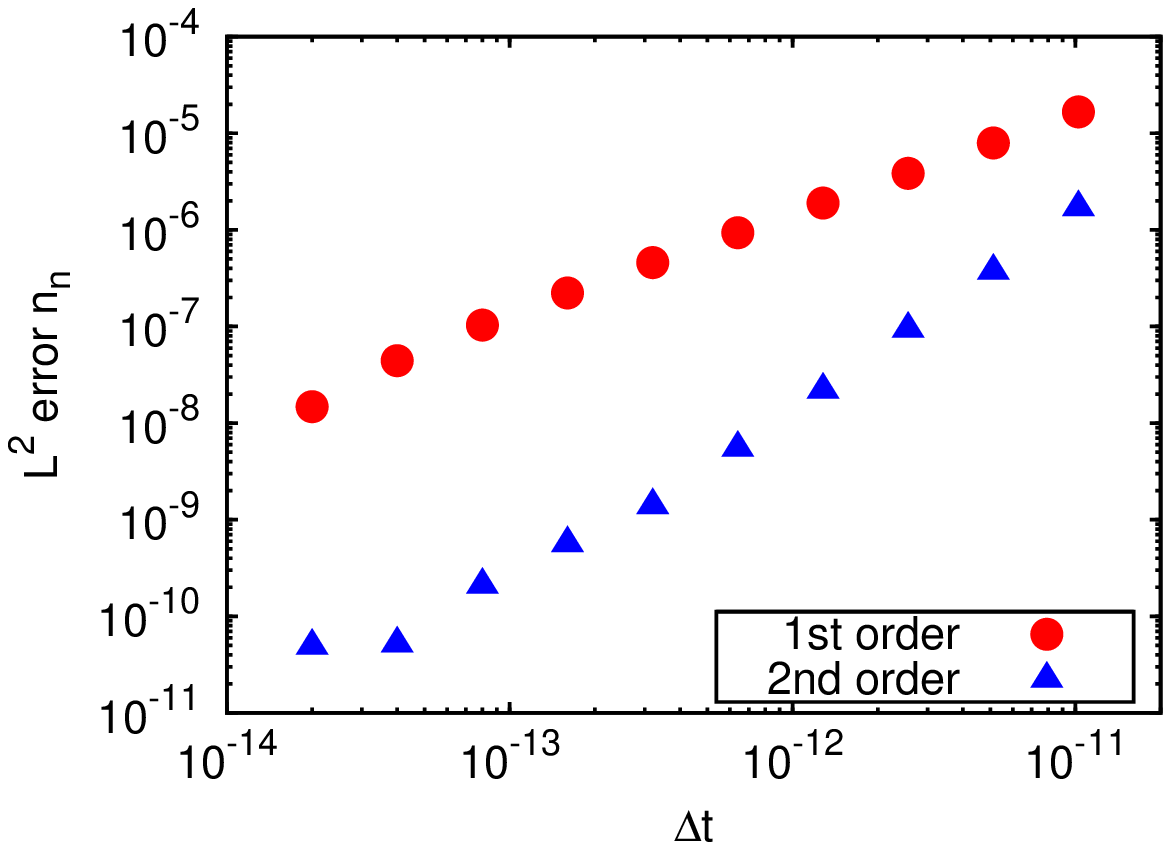}\hfill
   \caption{Normalized $L^2$ errors between the {\it reference} and the $\mathcal{T}_1$ (first order) and $\mathcal{T}_2$ (second order)
solutions for several decoupling time steps $\Delta t$ on a uniform grid of 4096 cells.
Top: electron (left) and positive ions (right); and bottom: negative ions.
             \label{order}}
 \end{center}
\end{figure}

First of all, we must verify the previous order estimates for the
$\mathcal{T}_1$ and $\mathcal{T}_2$ schemes given in Section~\ref{splitting_ef}.
We consider as initial condition the {\it reference} solution
at $t=10\,$ns.
In order to only evaluate errors coming from the 
decoupling techniques, $\mathcal{T}_1$ and $\mathcal{T}_2$,
we consider a fine splitting time step,
$\Delta t_{\rm s}=10^{-14}\,$s, to solve the drift-diffusion problem (\ref{trasp})
and a uniform grid;
then, we solve (\ref{semi_disc}) with both schemes
for several decoupling time steps $\Delta t_i$, and
calculate 
the normalized $L^2$ error between the first/second order
and {\it reference} solutions after time $t=2^{10}\Delta t_{\rm s}=1.024\times 10^{-11}\,$s.
Figure~\ref{order} shows results with 
$\Delta t_i=2^i\Delta t_{\rm s}$, where $i\in[1,10]$,
which clearly verify first and second order in time for
the $\mathcal{T}_1$ and $\mathcal{T}_2$ schemes, respectively,
and prove important gains in accuracy for same time steps.
For instance, for $\Delta t\leq10^{-12}\,$s 
the second order scheme provides solutions with $L^2$ errors
at least $100$ times lower than those obtained with the first order method.
\begin{figure}[!htb]
 \begin{center}
   \includegraphics[width=0.49\hsize]{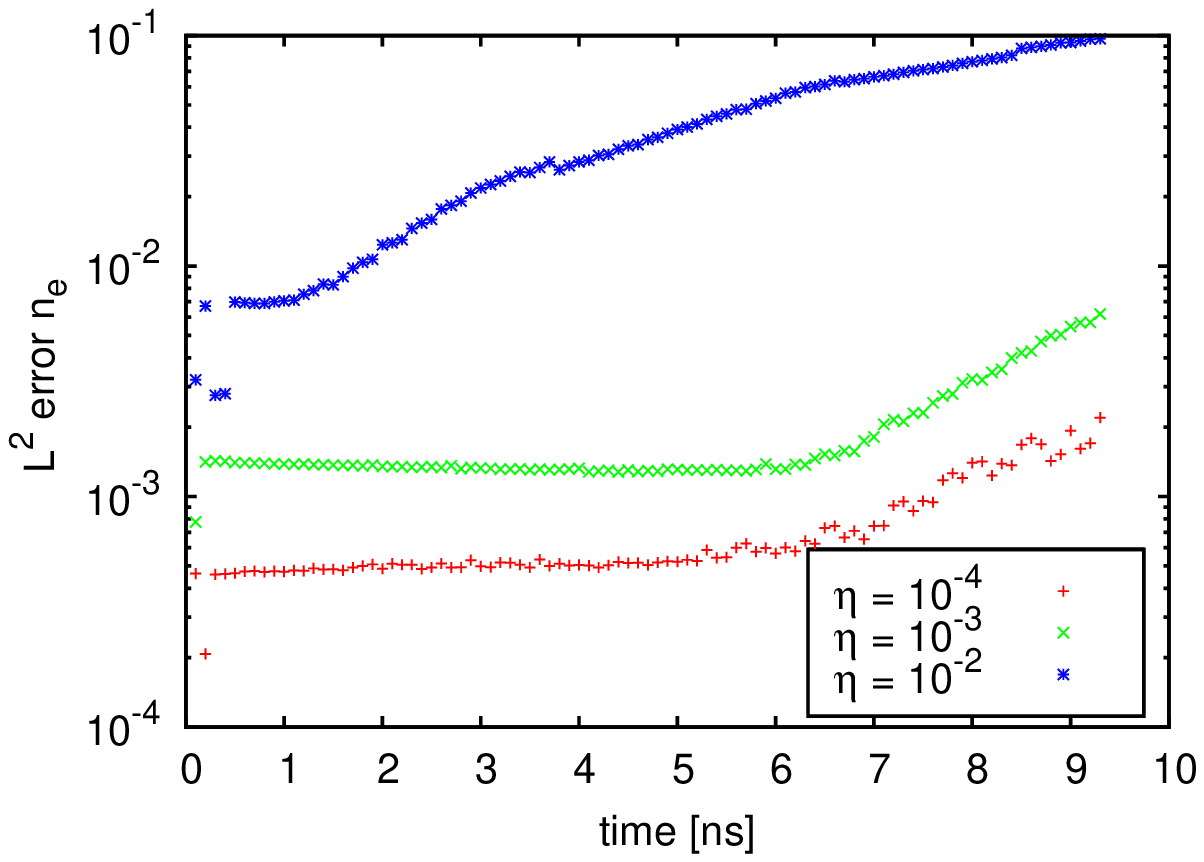}\hfill
  \includegraphics[width=0.49\hsize]{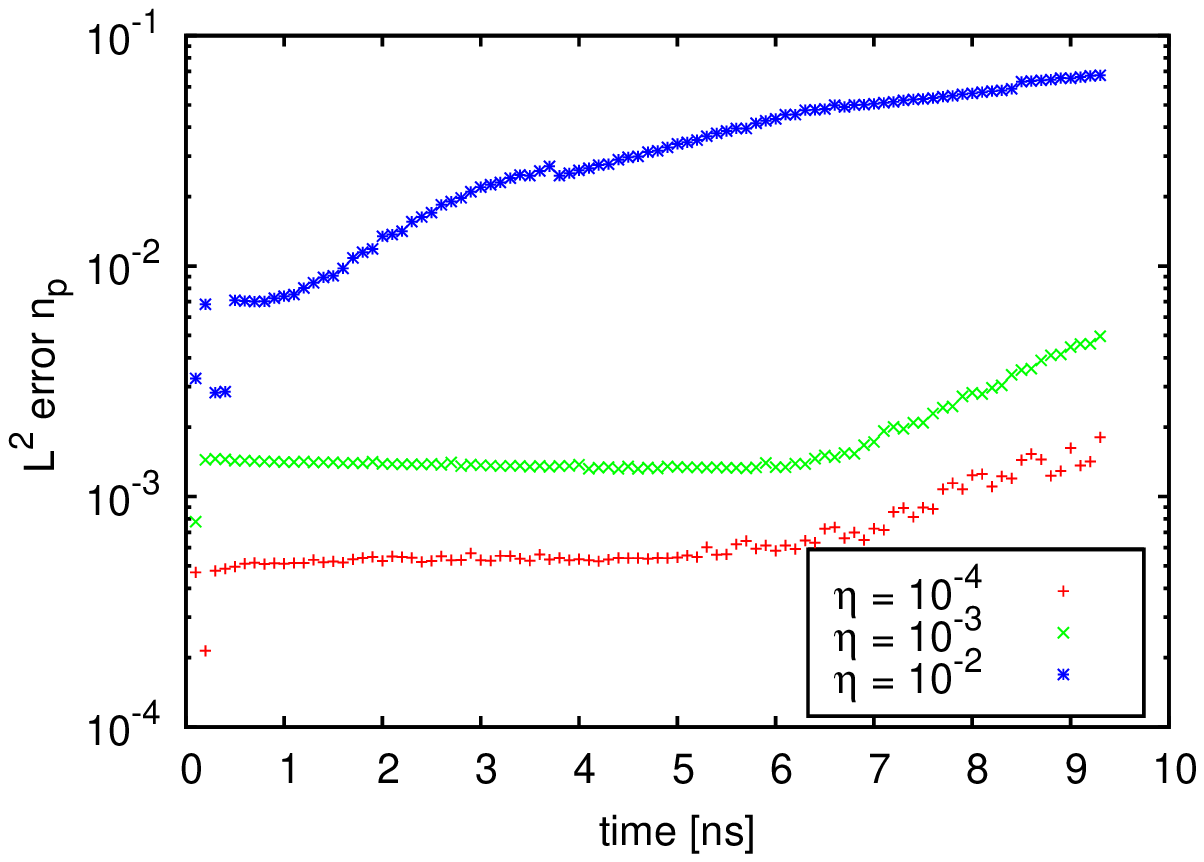}
 \includegraphics[width=0.49\hsize]{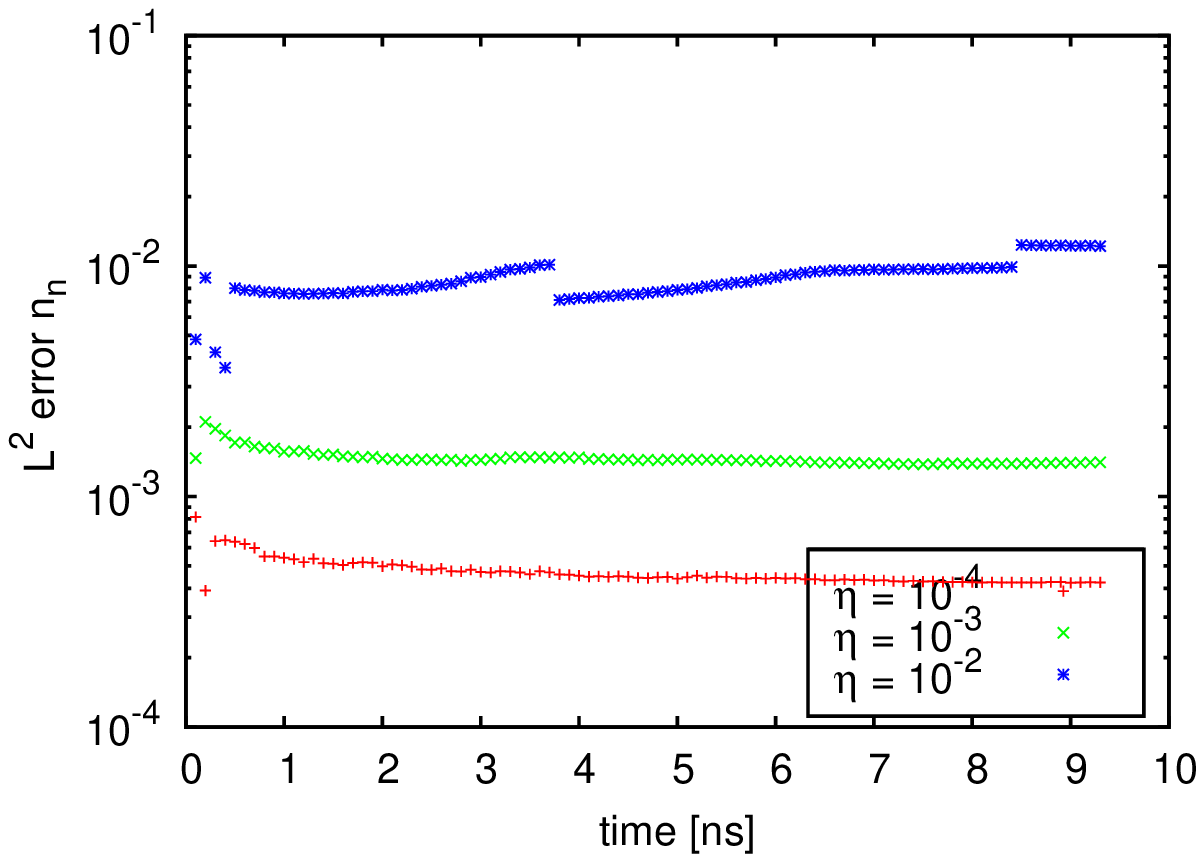}\hfill
   \caption{ Time evolution of the normalized $L^2$ errors between the {\it reference} and adapted solutions
with $\eta=\eta_{\mathcal{T}}=\eta_{\rm split}=\eta_{\rm MR}=10^{-4}$, $10^{-3}$, and $10^{-2},$
and 4096 cells corresponding to the finest discretization.
Top: electron (left) and positive ions (right); and bottom: negative ions.
 \label{error-12}}
 \end{center}
\end{figure}

Figure~\ref{error-12} shows the time evolution of the normalized $L^2$ error for each variable
between the time-space adapted
and {\it reference} solutions for several tolerances,
$\eta_{\mathcal{T}}=\eta_{\rm split}=\eta_{\rm MR}=10^{-4}$, $10^{-3}$, and $10^{-2}$.
These are rather approximations of the error since the {\it reference} and adapted solutions are not
evaluated exactly at the same time, and therefore, they are often slightly shifted
of about $\sim$$10^{-14}-10^{-13}$s.
In these tests, the decoupling time steps $\Delta t$ were limited by the dielectric relaxation time step,
$\Delta t_{DR}$,
after noticing an important amount of rejections of computed time steps
according to (\ref{delta_split_ef}), whenever $\Delta t \gtrsim$$1.5\times \Delta t_{DR}$.
Otherwise, $\Delta t$ is dynamically chosen in order to locally
satisfy the required accuracy $\eta_{\mathcal{T}}$, but it does not show important variations
considering the self-similar propagating phenomenon.
\begin{figure}[!htb]
 \begin{center}
   \includegraphics[width=0.49\hsize]{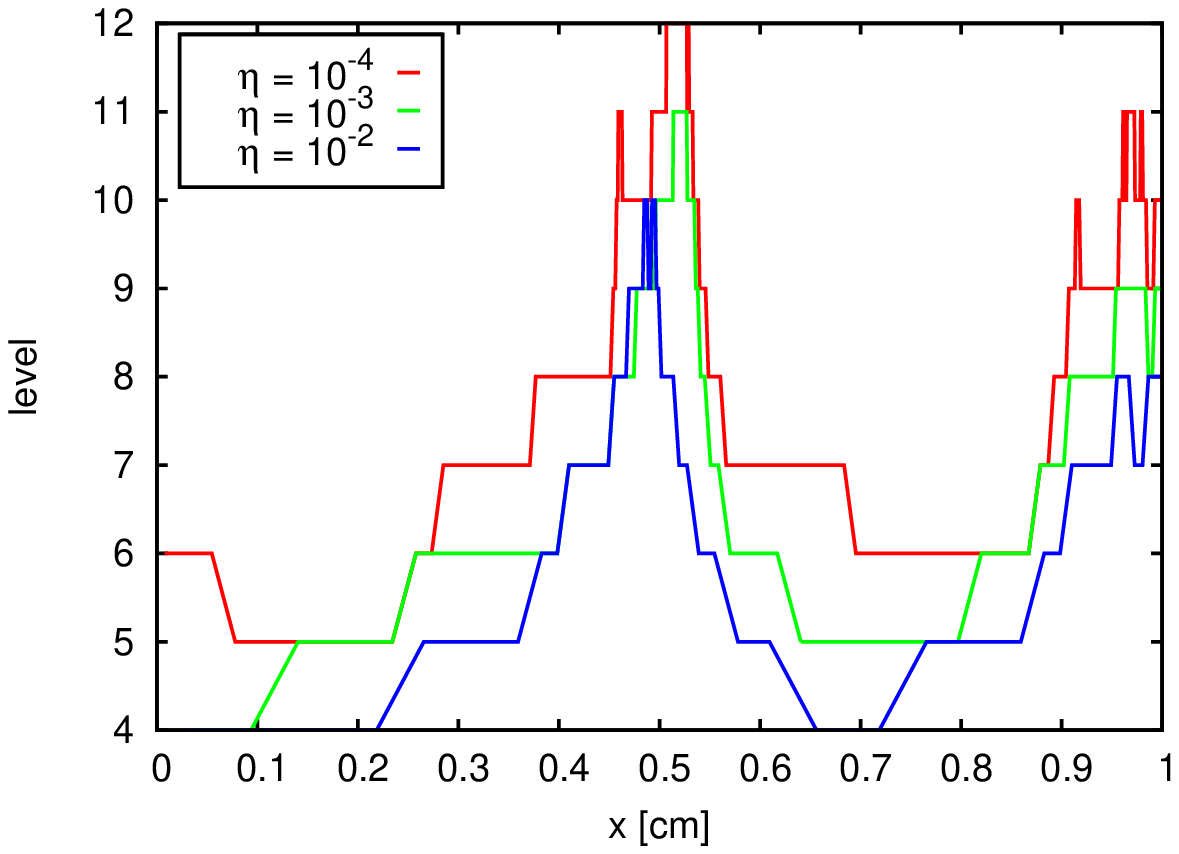}\hfill
 \includegraphics[width=0.49\hsize]{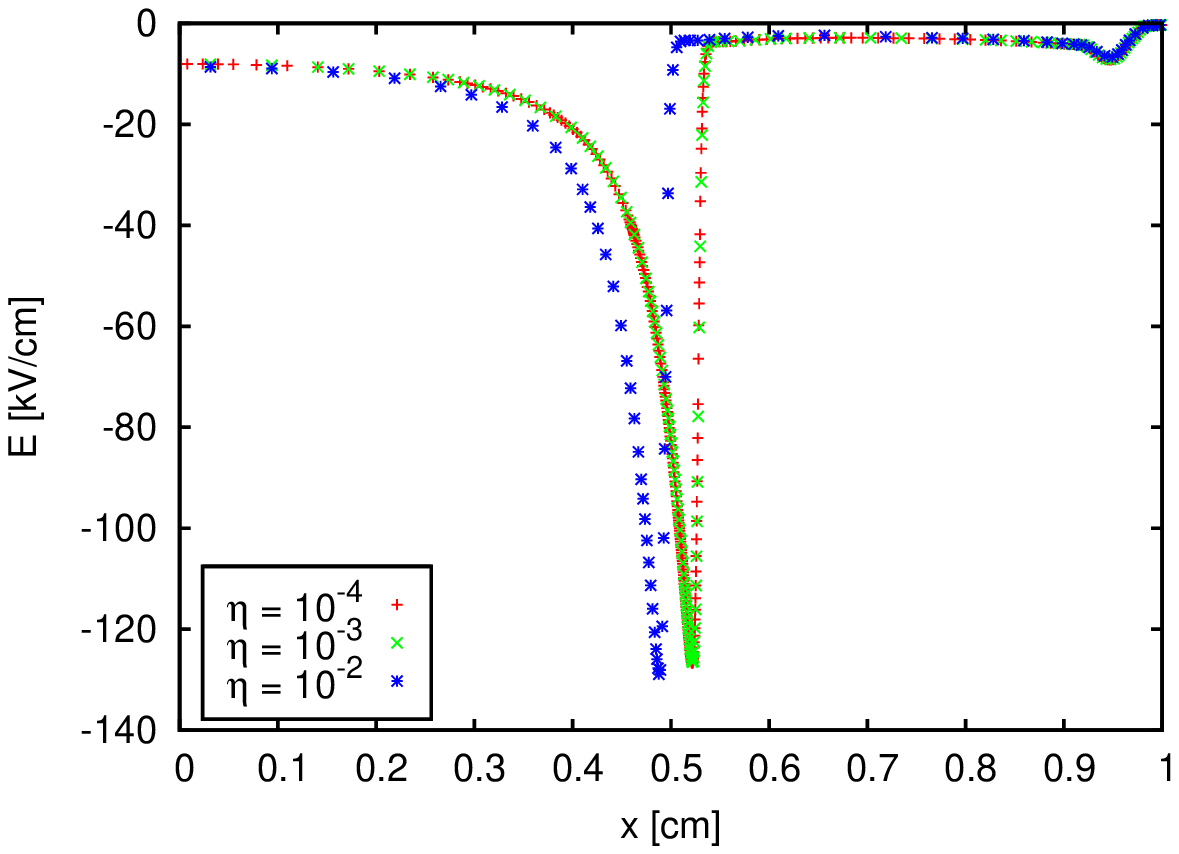}
 \includegraphics[width=0.49\hsize]{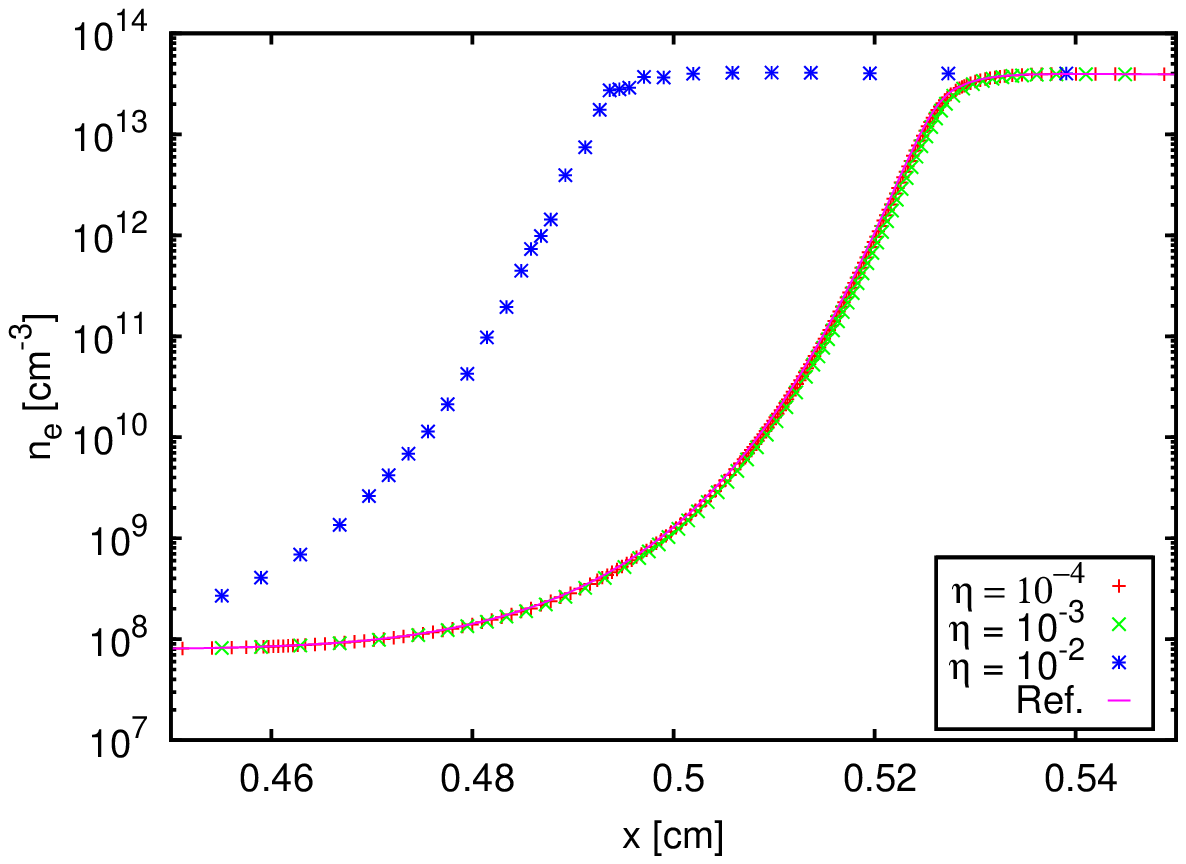}\hfill
 \includegraphics[width=0.49\hsize]{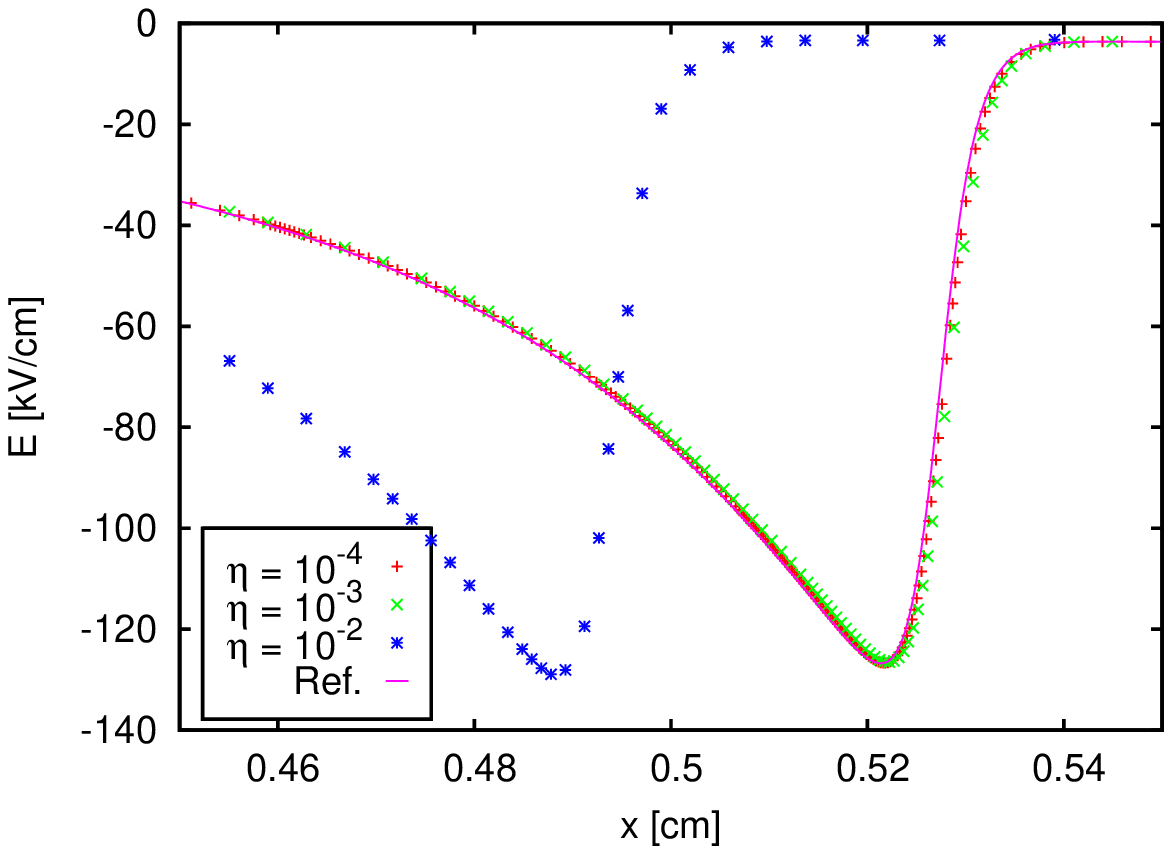}
   \caption{ Top: adapted grids (left)  and electric fields (right)
             at $t=8\,$ns with 4096 cells corresponding to the finest discretization, 
            and $\eta=\eta_{\mathcal{T}}=\eta_{\rm split}=\eta_{\rm MR}=10^{-4}$, $10^{-3}$, and $10^{-2}$.
Bottom: zoom on the electron distributions (left) and the electric field (right) with the same parameters,
and the {\it reference} solution.
 \label{level12}}
 \end{center}
\end{figure}

In Figure~\ref{level12}, we can see the corresponding adapted grid to
each previous configuration with different tolerances. 
The representation of the electric field and the densities
shows that for 
$\eta_{\mathcal{T}}=\eta_{\rm split}=\eta_{\rm MR}=10^{-2}$,
the streamer front propagates faster than in the {\it reference} case, with a slightly
higher peak of the electric field in the front.
On the other hand,
for $\eta_{\mathcal{T}}=\eta_{\rm split}=\eta_{\rm MR}\le10^{-3}$,
we observe a quite good agreement between
the adapted and {\it reference} resolutions.
\begin{figure}[!htb]
 \begin{center}
 \includegraphics[width=0.49\hsize]{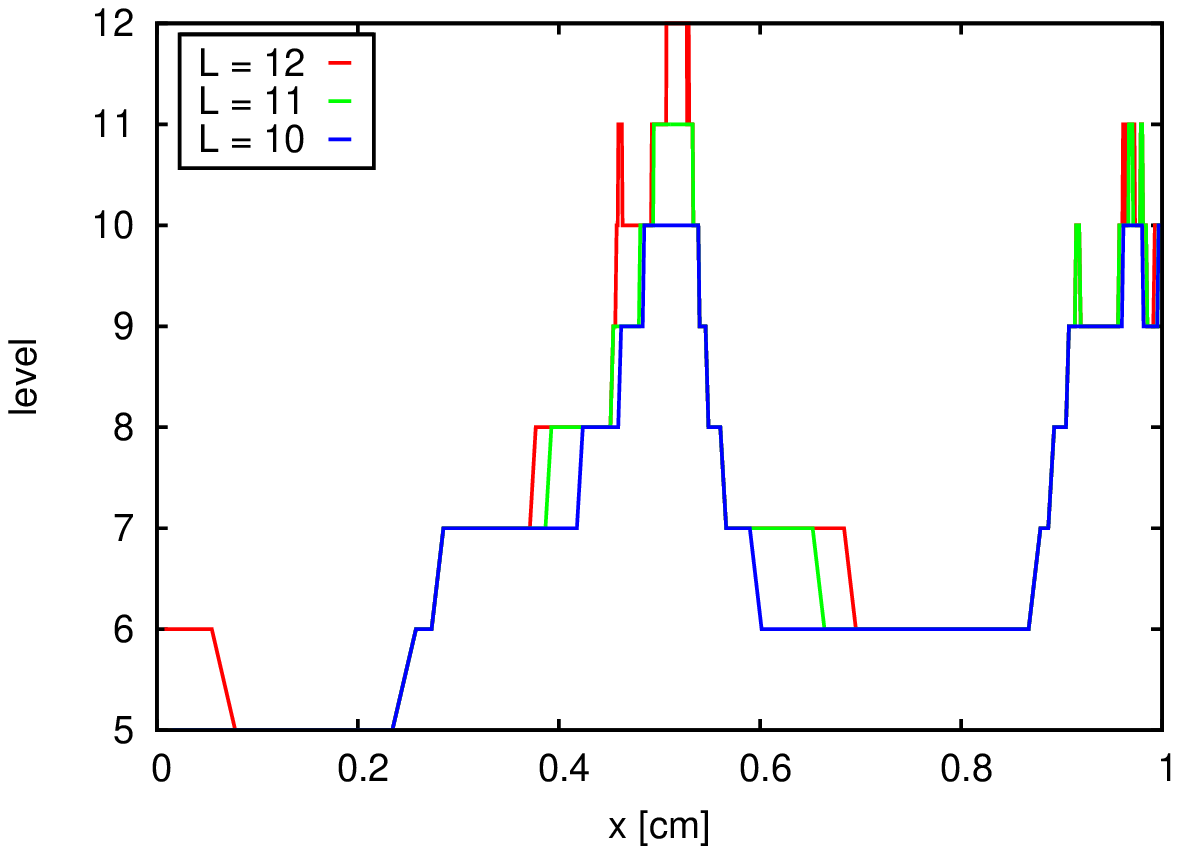}
 \includegraphics[width=0.49\hsize]{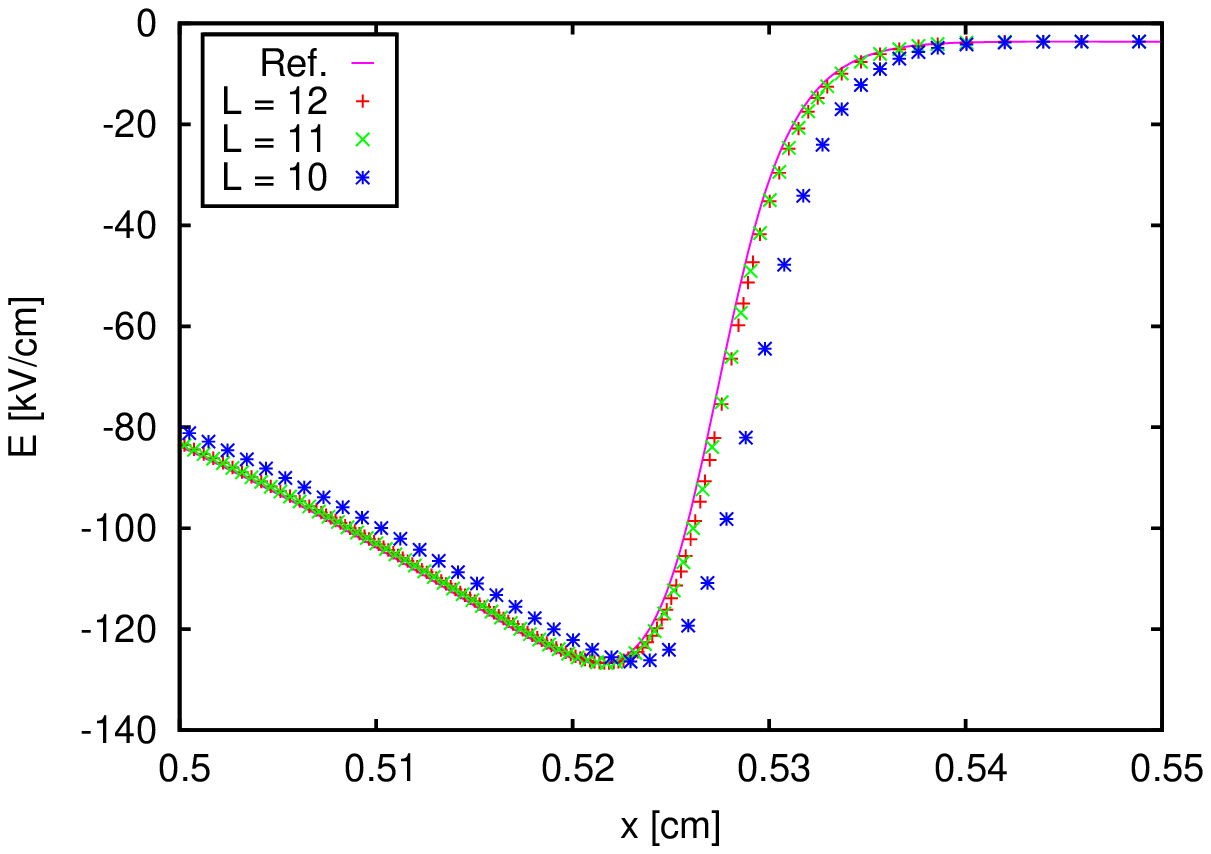}
   \caption{ Adapted grids (left)  and electric fields (right)
             at $t=8\,$ns, for several finest spatial discretization $L=10$, 11 and 12,
             $\eta_{\mathcal{T}}=\eta_{\rm split}=\eta_{\rm MR}=10^{-4}$, and the {\it reference} solution.
 \label{ef12-11-10}}
 \end{center}
\end{figure}

We consider now an accurate enough resolution with $\eta_{\mathcal{T}}=\eta_{\rm split}=\eta_{\rm MR}=10^{-4}$
and investigate the influence of the number of grids, that is, the 
finest spatial discretization at level $L$ that should be taken into account.
Figure~\ref{ef12-11-10} shows the adapted grids for 
$L=10$, $11$ and $12$, 
respectively equivalent to $1024$, $2048$ and $4096$ cells in the finest grid;
and a close-up of the corresponding
electric fields in the discharge head at $t=8\,$ns.
We see that for this level of tolerances,
the streamer front propagates slightly slower than the {\it reference} case
for $L=10$, whereas
$L=11$ gives
already good resolutions compared with the {\it reference} solution and
with $L=12$.
In particular, higher values of $L$ would need lower tolerances in order
to retain regions at the finest level;
this is already the case for $L=13$
(equivalent to $8192$ cells).
Therefore, $L=11$  with $2048$ cells at the finest level
seems to be an appropriate choice
for this level of accuracy.

Table~\ref{Table11} 
summarizes the number of cells in the adapted grid ($\#AG$) at time  $t=8$ns,
and the corresponding data compression ($DC$) defined as the percentage of active cells
with respect to the equivalent number of cells for the finest discretization, 
in this case $2048$ for $L=11$.
For this propagating case, the data compression remains of the same order during
the time simulation interval.
The CPU computing times correspond to a time domain of study of
$t\in[0,10]\,$ns computed by one sole CPU core.
If we consider for example
total computing time for $L=11$ and 
tolerances $\eta_{\mathcal{T}}=\eta_{\rm split}=\eta_{\rm MR}=10^{-4}$,
it is $\sim$$44$ times less expensive 
with respect to a resolution on
a uniform grid with $2048$ cells and $\eta_{\mathcal{T}}=\eta_{\rm split}=10^{-4}$
(CPU time of $8552\,$s).
This is quite reasonable, taking into account that
the computing time for the electric
field resolution is proportional to 
at least $\mathcal{O}(N^2)$ 
for $N$ computing cells, after
(\ref{EFi}).
 \begin{table}[!htb]
\caption{Number of cells in the adapted grid ($\#AG$) and data compression ($DC$) 
at time  $t=8$ns, CPU computing time for $t\in[0,10]$ ns, $L=11$, and 
several tolerances $\eta=\eta_{\mathcal{T}}=\eta_{\rm split}=\eta_{\rm MR}$. \label{Table11}}
\smallskip
\centerline{%
\bgroup\offinterlineskip
\vbox{%
\halign to 0.45\hsize{%
\vrule height 12pt  width 0pt%
\kern 11pt# &#\kern12pt&
\hfill# \hfill&#\kern12pt&
\hfill# \hfill&#\kern12pt&
#\hfil\cr
\noalign{\smallskip\hrule}
  $\eta$ &&  $\#AG$  && $DC$\% && CPU(s)\cr
\noalign{\smallskip\hrule}
  $10^{-6}$  &&   $724$           &&  $35.35$     && $1360$\cr
  $10^{-5}$  &&   $421$           &&  $20.56$     && $\phantom{0}517$\cr
  $10^{-4}$  &&   $263$           &&  $12.84$     && $\phantom{0}193$\cr
  $10^{-3}$  &&   $138$           &&  $\phantom{0}6.74$          && $\phantom{00}66$\cr
  $10^{-2}$  &&   $\phantom{0}70$            &&  $\phantom{0}3.42 $          && $\phantom{00}24$\cr
\noalign{\smallskip\hrule}}
}\egroup
}
\end{table}

In conclusion, in this section we have shown that the numerical strategy developed can be efficiently 
applied to simulate the propagation of highly nonlinear ionizing waves as streamer discharges. 
An important reduction of computing time results from significant data compression
with still accurate resolutions.
In addition, 
this study allows to properly tune the various simulation parameters in order
to guarantee a fine resolution of more complex configurations,
based on the time-space accuracy control 
capabilities of the method.

\FloatBarrier

\subsection{Simulation of multi-pulsed discharges}
In this section, we analyze the performance of the proposed numerical strategy on the simulation
of nanosecond repetitively pulsed discharges \cite{Pilla:2006,pai:2010}. 
The applied voltage profile for this
type of discharges is a high voltage {\it pulse} followed by a zero voltage {\it relaxation}
phase. The typical pulse duration is $\sim$$10^{-8}\,$s, while the relaxation phase
takes over $\sim$$10^{-4}\,$s. 
The detailed experimental study of these discharges in air
has shown that the cumulative effect of repeated pulsing achieves a steady-state
behavior \cite{pai:2010}. 
In the following illustrations,
we choose a pulse duration of $T_{\rm p}=15\,$ns, which is  approximately
equal to the time that is needed for the discharge to cross the inter-electrode gap.
The rise time considers the time needed to go from zero to the maximum voltage and
it is set to $T_{\rm r}=2\,$ns.
The pulse repetition period is set to $T_{\rm P} = 10^{-4}\,$s, equal to
$10\,$kHz of repetition frequency, a typical value
used in experiments \cite{Pilla:2006}.
We model the voltage pulse $P$ by using sigmoid functions
\begin{equation}\label{sigmoid}
P(t,s,r,p)= 1 - \sigma(-t,-s,r) - \sigma(t,s+p,r),
\end{equation}
with
\begin{equation}\label{sigma}
\sigma(t,s,r)=\frac{1}{1+\exp(-8(t-s)/r)},
\end{equation}
for time $t$, where $s$ indicates when the pulse starts; 
$r$ is the rise time; and $p$ is the pulse duration;
$t$, $s$, $r$, $p \in [0,T_{\rm P}]$. 
With a maximum applied voltage $V_{\rm max}$,
the applied voltage $V(t)$ is computed by
\begin{equation}
V(t)=V_{\rm max}\cdot P \left( t-\left\lfloor \frac{t}{T_{\rm P}} \right\rfloor \cdot T_{\rm P}, T_{\rm r},T_{\rm r},T_{\rm p}\right).
\end{equation}
In repetitively pulse discharges at atmospheric pressure and $300\,$K,
as discussed in \cite{Pancheshnyi:2005b,Wormeester:2010},
electrons attach rapidly to $O_2$ molecules during the interpulse
to form negative ions (characteristic time scale of $20\,$ns).
Then, the rate of the plasma decay is determined by ion-ion
recombination \cite{Kossyi:1992,Pancheshnyi:2005b,Wormeester:2010}.
When the next voltage pulse is applied, electrons
are detached with a rate 
taken from \cite{Benilov:2003}.
Therefore, as initial condition we assume a distribution similar
to the end of the interpulse phase with a 
homogeneous preionization consisting of positive and negative ions with a density of $10^{9}\,$cm$^{-3}$.
For electrons, we consider a low homogeneous background of $10^{1}\,$cm$^{-3}$.
This small amount of electrons as initial condition has a negligible influence on the results.

\begin{figure}[!htb]
 \begin{center}
 \includegraphics[width=0.49\hsize]{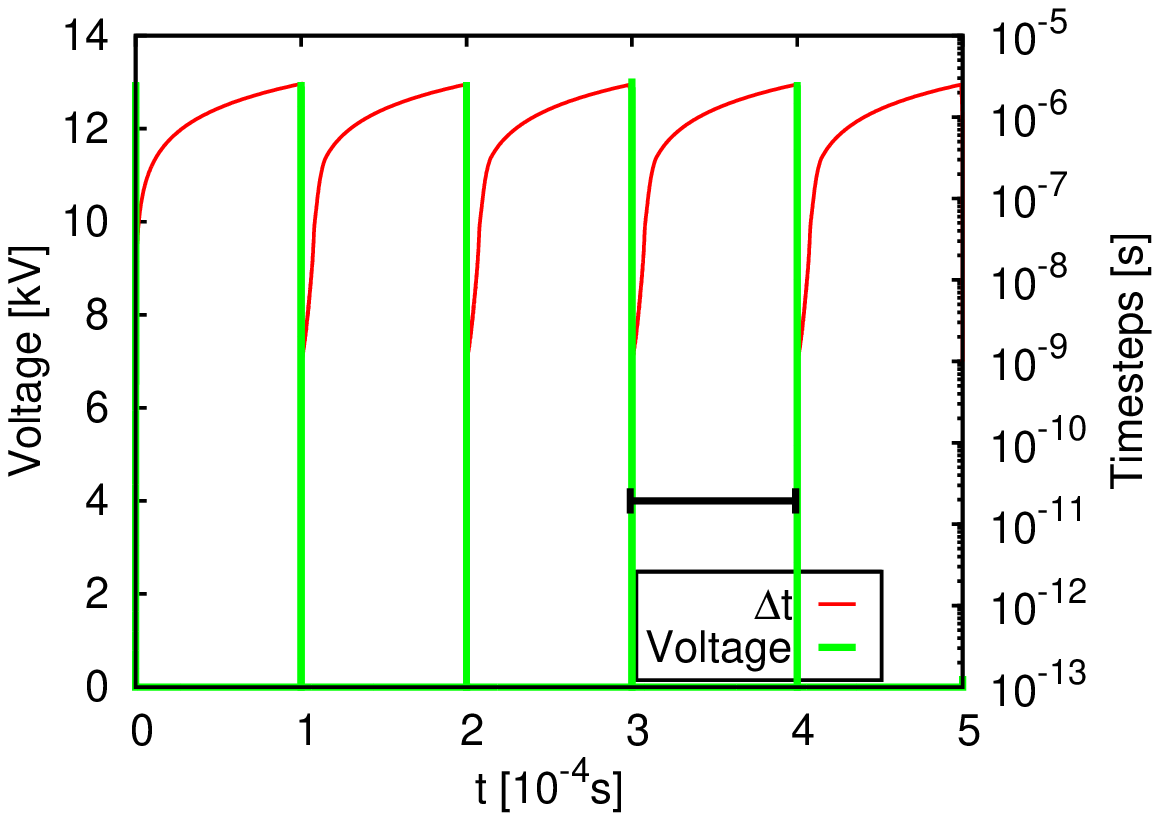}
 \includegraphics[width=0.49\hsize]{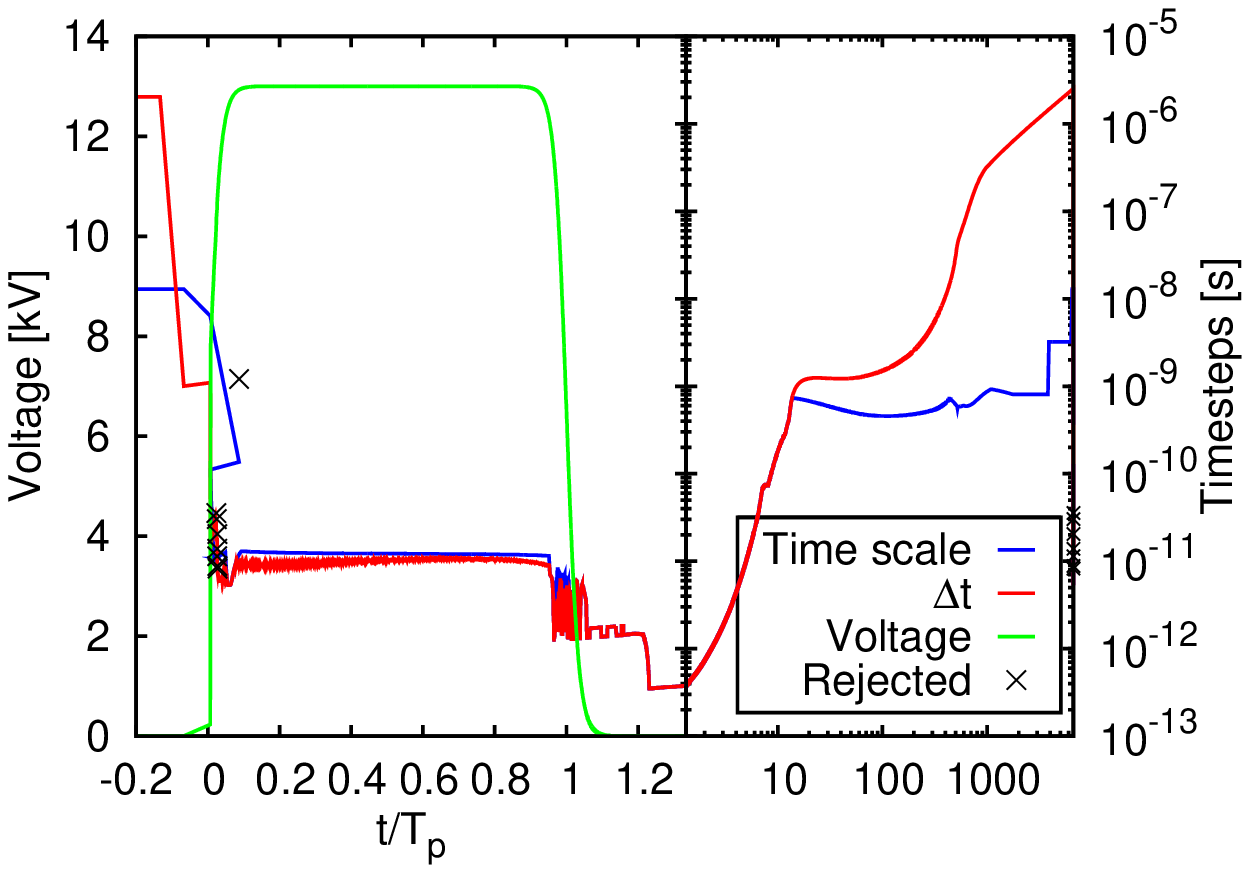}
   \caption{Time evolution of the applied voltage and the decoupling time steps $\Delta t$ 
for a multi-pulse simulation
for the first 6 pulses (left) and for the 4th one (right) with its subsequent relaxation.
            Rejected time steps are marked with black crosses, while the minimum time scale corresponds to the blue line.
   \label{multipulse-timestep}}
 \end{center}
\end{figure}

We set the tolerances to $\eta_{\mathcal{T}}=\eta_{\rm split}=\eta_{\rm MR}=10^{-4}$
and consider $L=11$ grid levels,
equivalent to $2048$ cells in the finest grid. 
As in the previous configuration, homogeneous Neumann boundary conditions were considered for the
drift-diffusion equations.
Figure~\ref{multipulse-timestep} shows the time evolution of the decoupling time steps
and the applied voltage for the first six pulses, even though
simulation was performed for 100 pulses, that is
$t\in[0,10^{-2}]\,$s.
This simulation took over $8$h$44$m 
while running the electric field
computation in parallel on $6$ CPU cores
of the same AMD Opteron 6136 Processor cluster;
this gives an average of $5.24$ minutes per pulse period.
Figure~\ref{multipulse-timestep} shows also the fourth pulse for which
the steady-state of the periodic phenomenon
was already reached and almost the same numerical performance is
reproduced during the rest of computations.
The time steps are about $\sim$$10^{-11}\,$s during pulses,
then increase from 
$\sim$$10^{-12}\,$s
up to about $\sim$$10^{-6}\,$s during a period
$\sim$$6000$ times longer,
for which standard stability constraints are widely overcome
according to the required accuracy tolerance.
Solving this problem for such different scales with a constant
time step is out of question and even a standard strategy
that considers the minimum of all time scales would limit considerably the efficiency of
the method as it is shown in the representation.
In this particular case, the dielectric relaxation is the governing time scale 
during the discharge as in the previous case with constant applied voltage,
whereas the post-discharge phase is alternatively ruled by diffusive or convective CFL,
or by ionization time scale, with all security factors and
CFL conditions set to one in Figure~\ref{multipulse-timestep}.

The computation is initialized with a time step included in the pulse
duration.
Nevertheless, after each relaxation phase, since
the new time step is computed based on the previous one
according to (\ref{delta_split_ef}), this new time step
will surely skip the next pulse.
In order to avoid this, each time we get into a new
period, that is when $\lfloor t/T_{\rm P}\rfloor$ changes, we
initialize the time step with $\Delta t = 0.5T_{\rm r}=1\,$ns.
This time step is obviously rejected as seen in 
Figure~\ref{multipulse-timestep}, as well as the next ones, until
we are able to retrieve the right dynamics of the phenomenon
for the required accuracy tolerance.
No other intervention is needed neither for modeling parameters
nor for numerical solvers in order to automatically
adapt the time step needed to describe the various time scales of
the phenomenon within a prescribed accuracy.

\begin{figure}[!htb]
 \begin{center}
 \includegraphics[width=0.49\hsize]{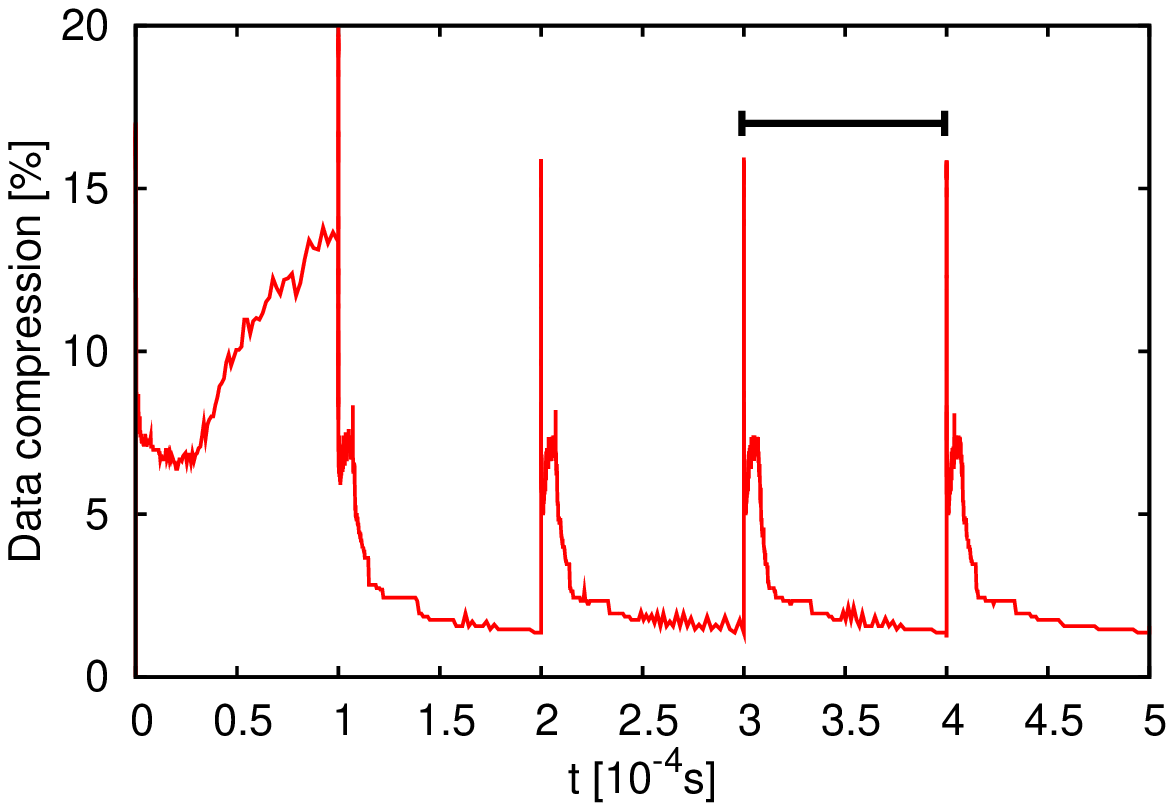}
 \includegraphics[width=0.49\hsize]{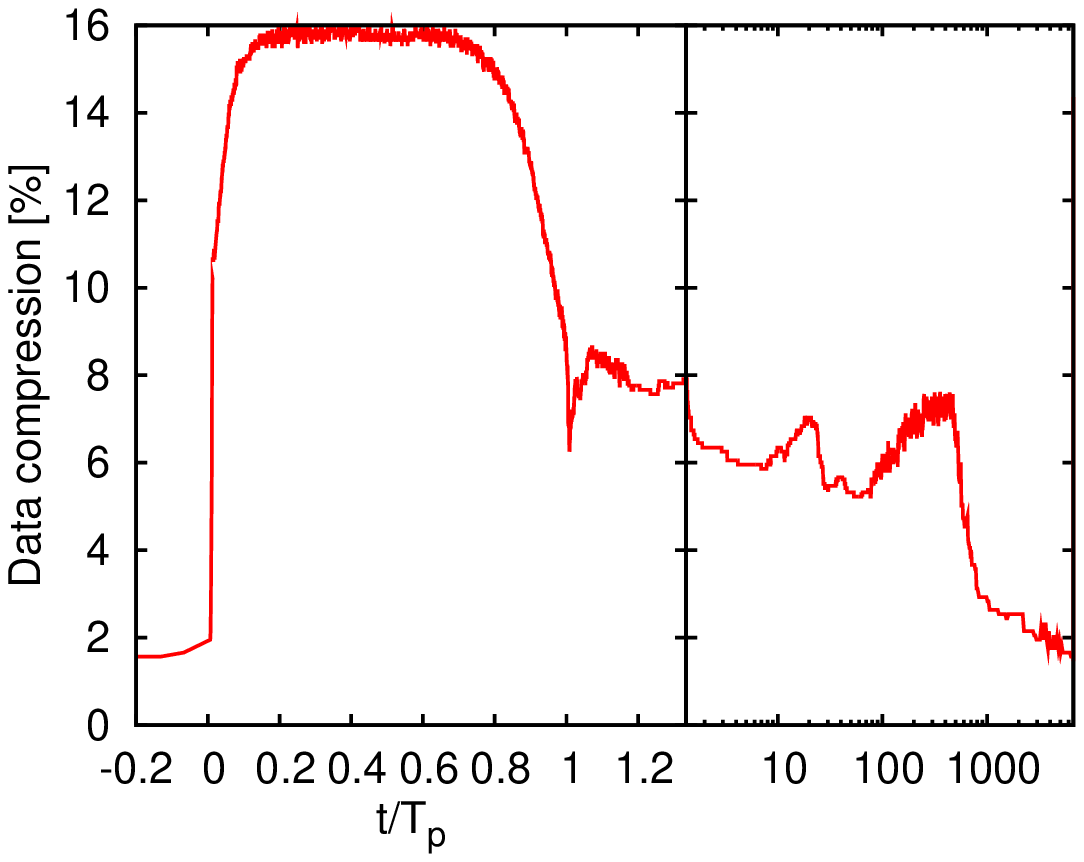}
   \caption{Time evolution of the data compression for a multi-pulse simulation
for the first six pulses (left) and for the fourth one (right) with its subsequent relaxation.
   \label{multipulse-compression}}
 \end{center}
\end{figure}

Figure~\ref{multipulse-compression} represents the time evolution of the data 
compression which ranges from $\sim$$2\%$ up to $\sim$$16\%$
during each pulse period.
Regarding only the electric field resolution with the same time
integration strategy, a grid adaptation technique involves resolutions
$\sim$$39$ to $\sim$$2500$ times faster,
based on a really rough estimate for $\mathcal{O}(N^2)$ operations.

\begin{figure}[!ht]
 \begin{center}
  \includegraphics[width=0.44\hsize]{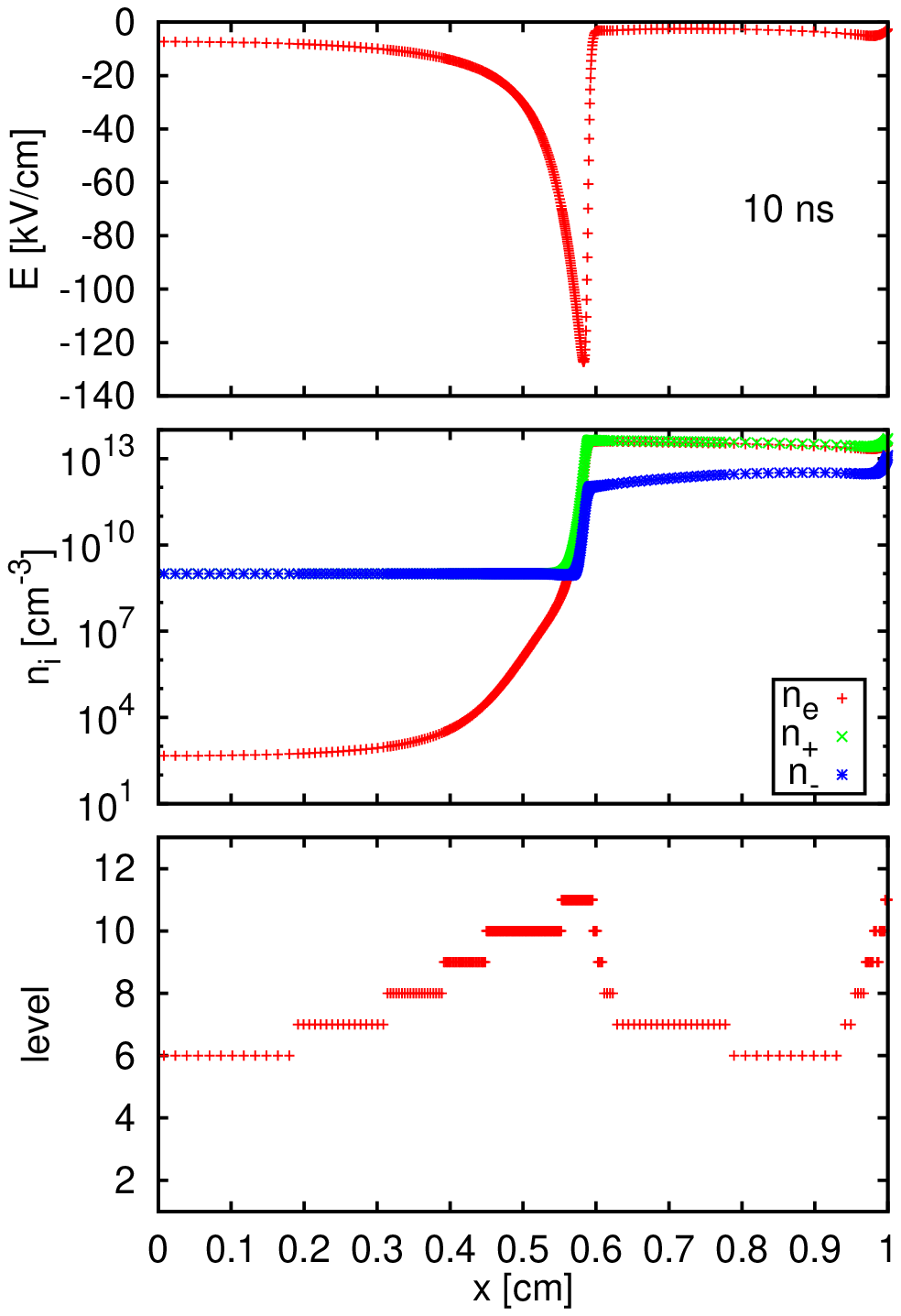}
 \includegraphics[width=0.44\hsize]{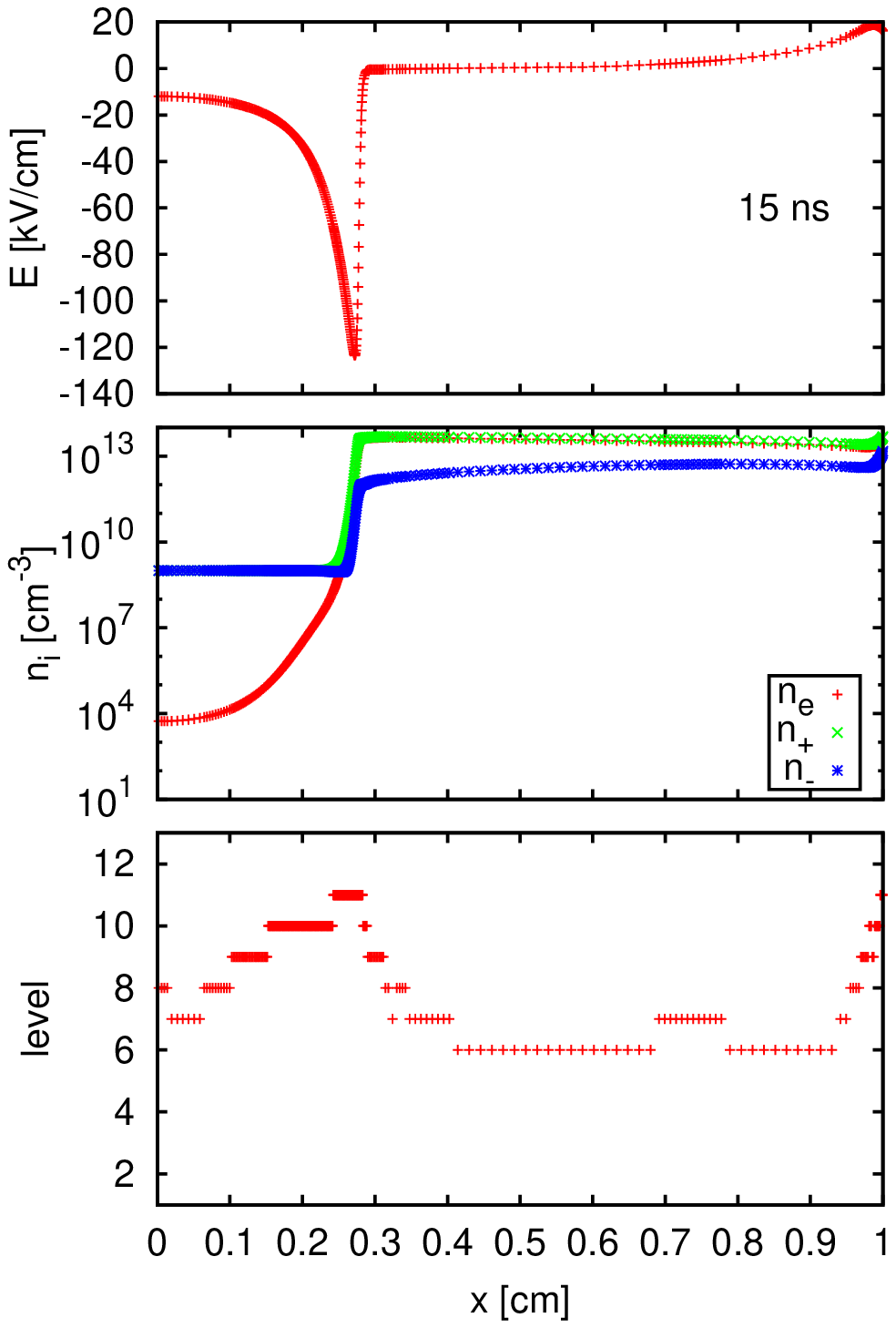} 
   \includegraphics[width=0.44\hsize]{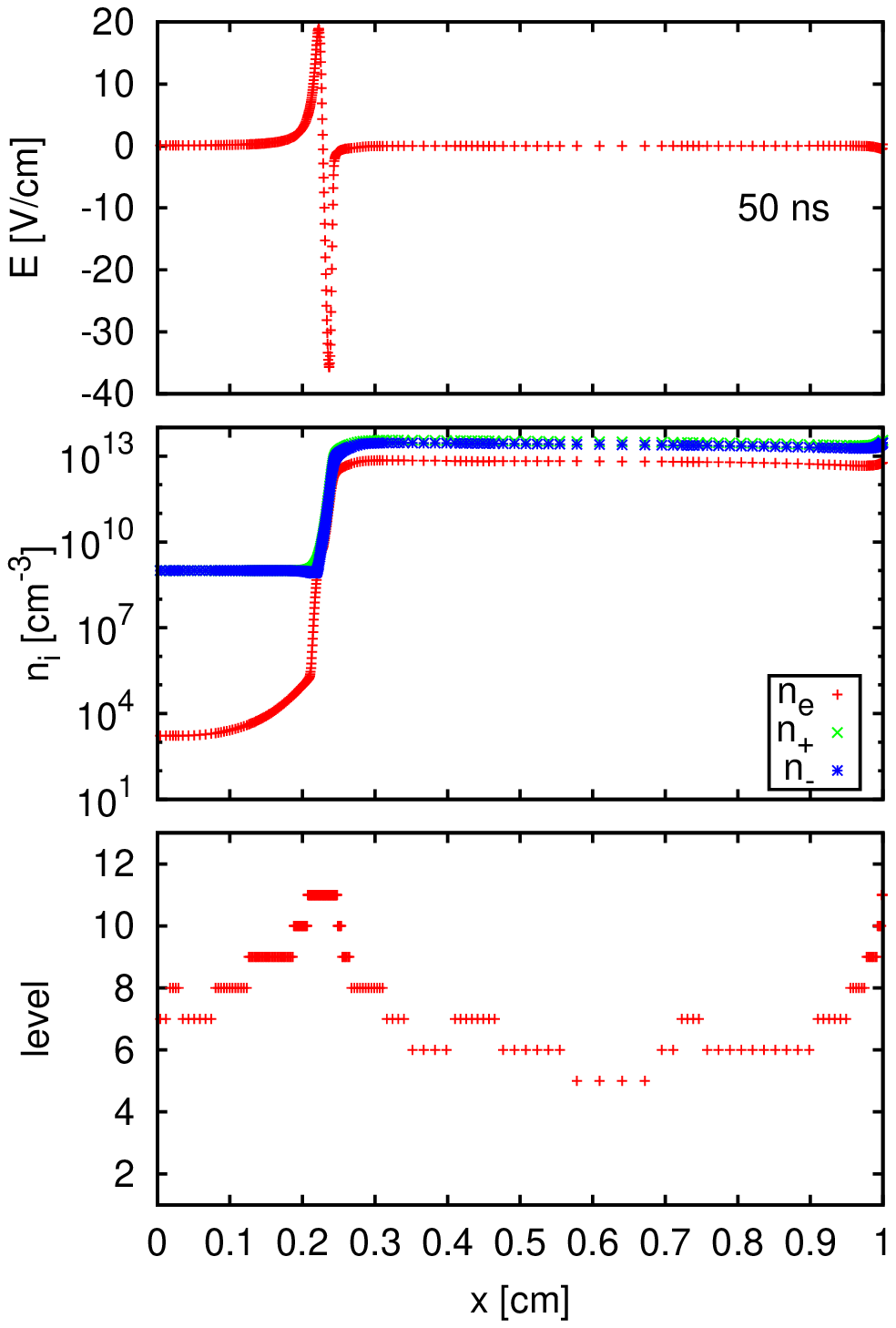}
\includegraphics[width=0.44\hsize]{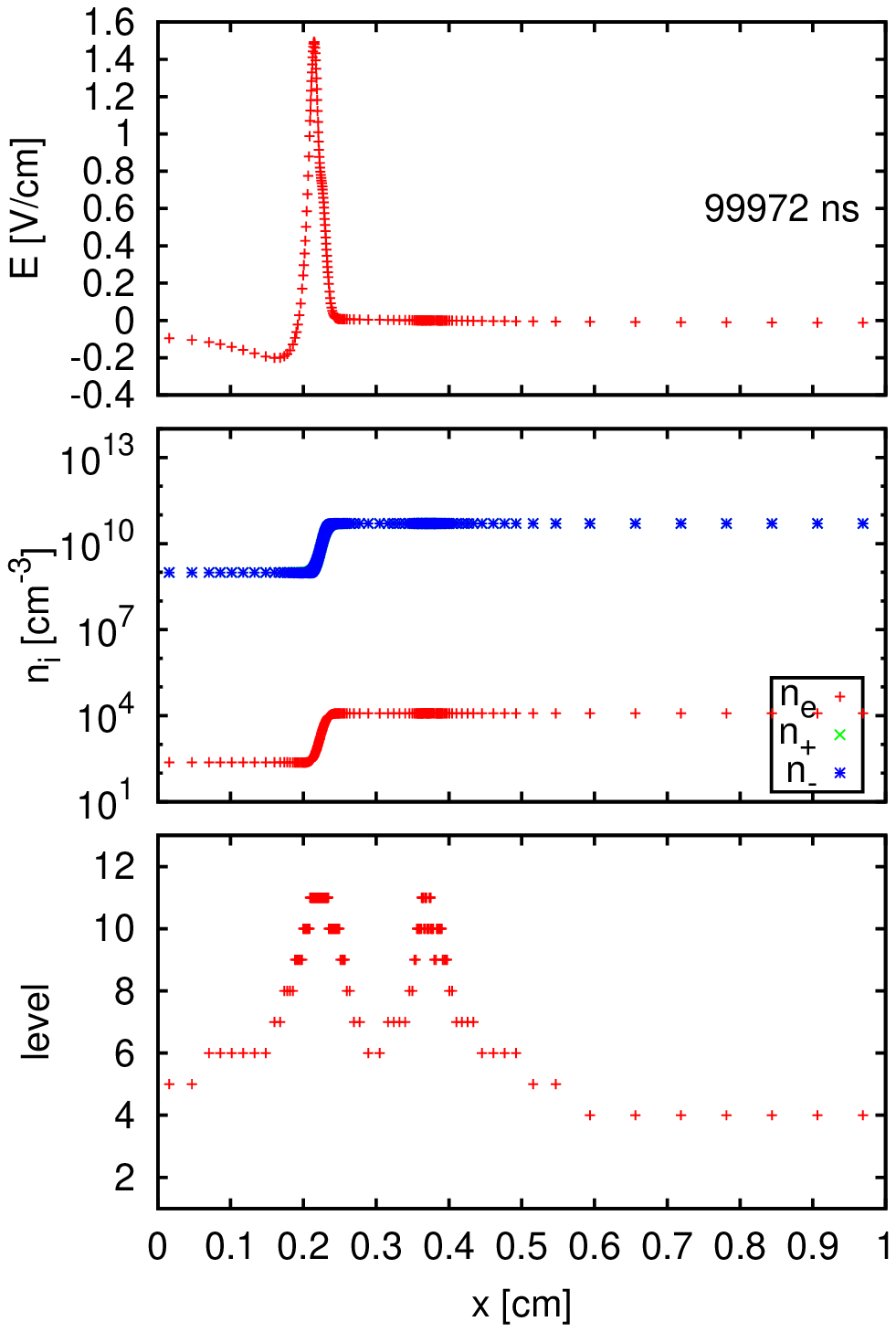}
   \caption{First period of pulsed discharges.
            Top: propagation of the discharge in the domain at $t=10\,$ns after the beginning of the pulse (left);
            and at $t=15\,$ns (right). 
            Bottom: relaxation on the short time scale $t=50\,$ns;
            and end of the relaxation phase after $t=99972\,$ns (right).
 \label{multipulse-domain-1}}
 \end{center}
\end{figure}

Figure~\ref{multipulse-domain-1} presents the discharge dynamics for the first period. 
First, we observe at $t=10\,$ns after the beginning of the pulse, the propagation of a positive streamer in the gap. 
In Section~\ref{sec-exact}, a preionization of positive ions and electrons was used to ensure the positive streamer propagation. 
In this section, seed electrons ahead of the streamer front are created as the front propagates by detachment of 
negative ions initially present. We note that at $15\,$ns, which corresponds to almost the end of the plateau before 
the decrease of the applied voltage, the discharge has crossed $\sim$$0.75\,$cm of the $1\,$cm gap. 
As a consequence, during the voltage decrease and at the beginning of the relaxation phase where the applied voltage is zero, 
there is a remaining space charge and steep gradients of charged species densities in the gap. 
Then for $t=50\,$ns, Figure~\ref{multipulse-domain-1} shows that the electric field in the discharge is almost equal 
to zero except in a small area where steep gradients of the electric field are observed
but with peak values of only $30\,$V/cm.
We have checked that this area corresponds to the location of the streamer head at the end of its propagation
as it is seen in the representation. 
We note that in the post-discharge, electrons are attaching and then at $t=50\,$ns, 
the density of positive ions is almost equal to the density of negative ions in the whole gap. 
At $t=99972\,$ns, the densities of charged species have significantly decreased due to charged species recombination. 
However, it is interesting to note that the location of the previous streamer head can still be observed at the same location 
as at $t=50\,$ns, but with much smaller gradients of charged species densities and a very small electric field. 
This final state is the initial condition of the second pulse with a non-uniform axial preionization with positive and negative ions
and a much smaller density of electrons.
\begin{figure}[!ht]
 \begin{center}
   \includegraphics[width=0.44\hsize]{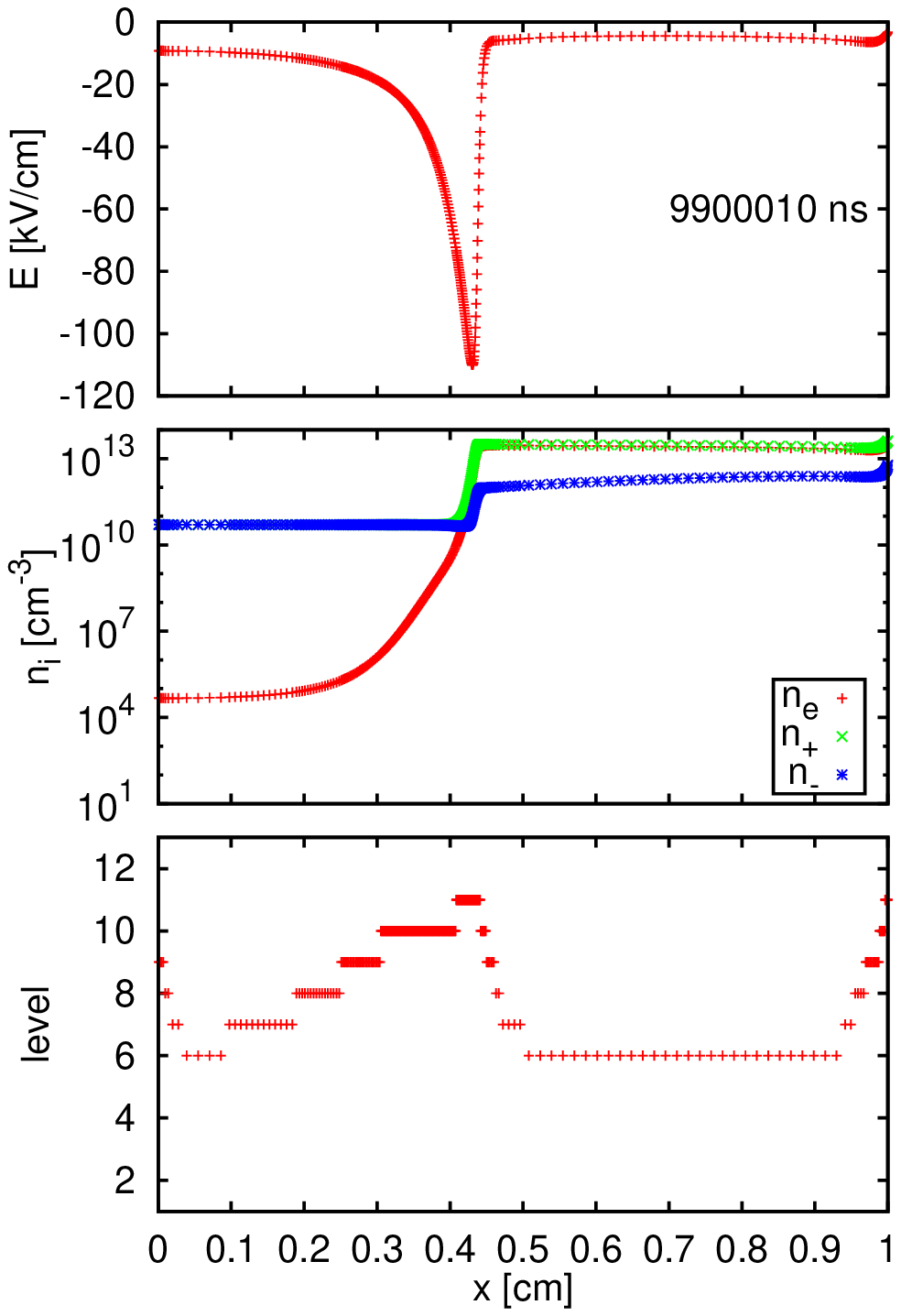}
 \includegraphics[width=0.44\hsize]{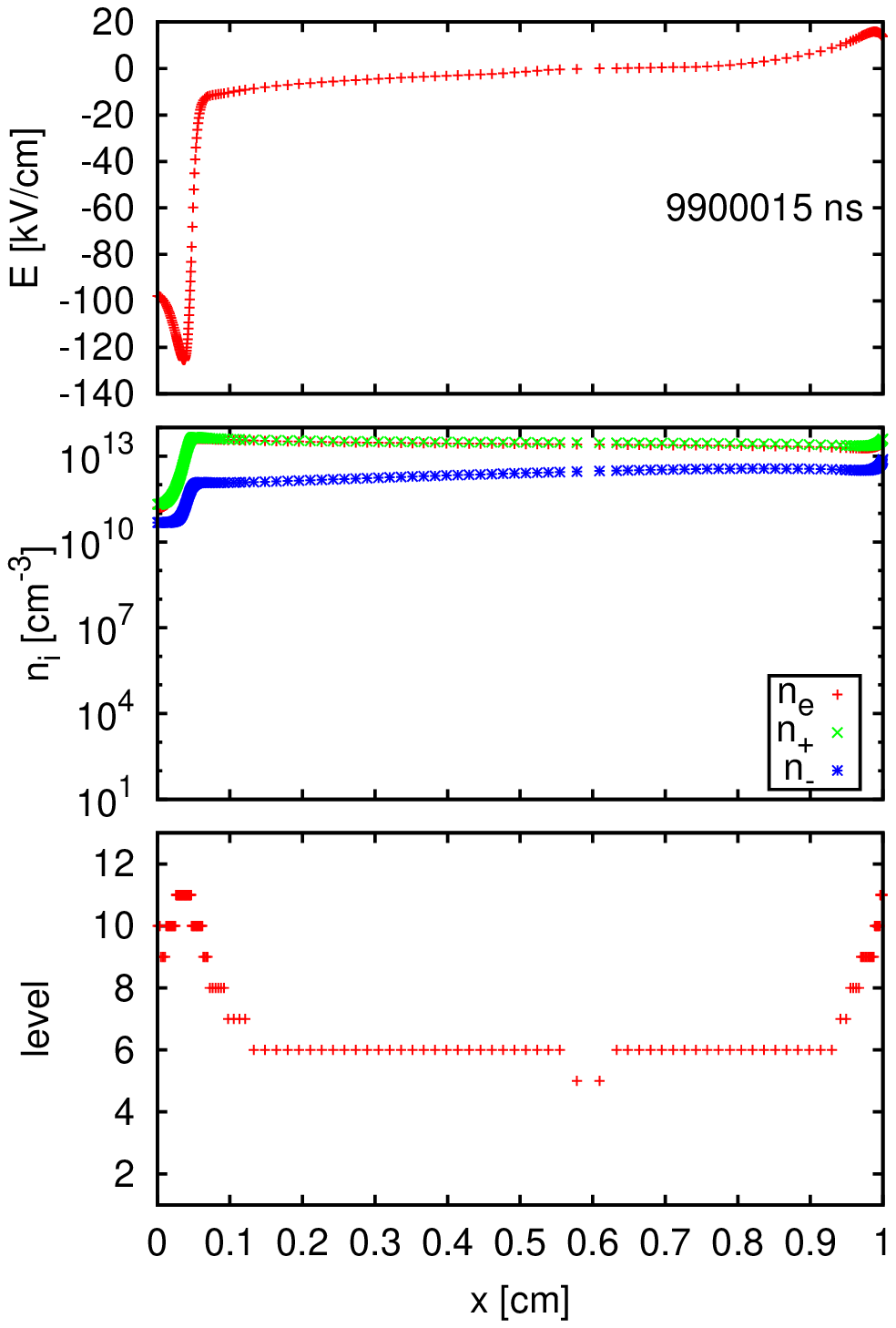}  
 \includegraphics[width=0.44\hsize]{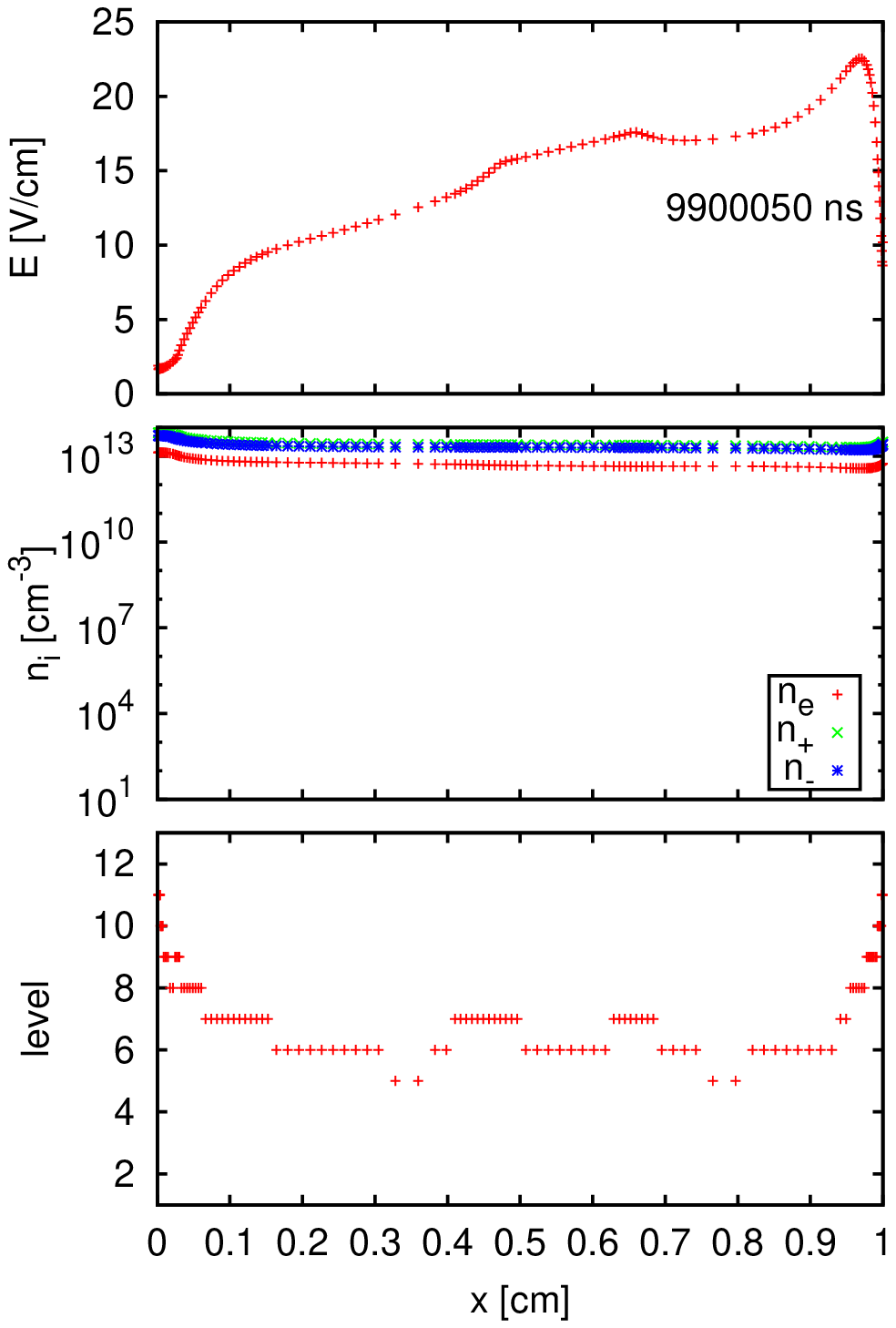}
 \includegraphics[width=0.44\hsize]{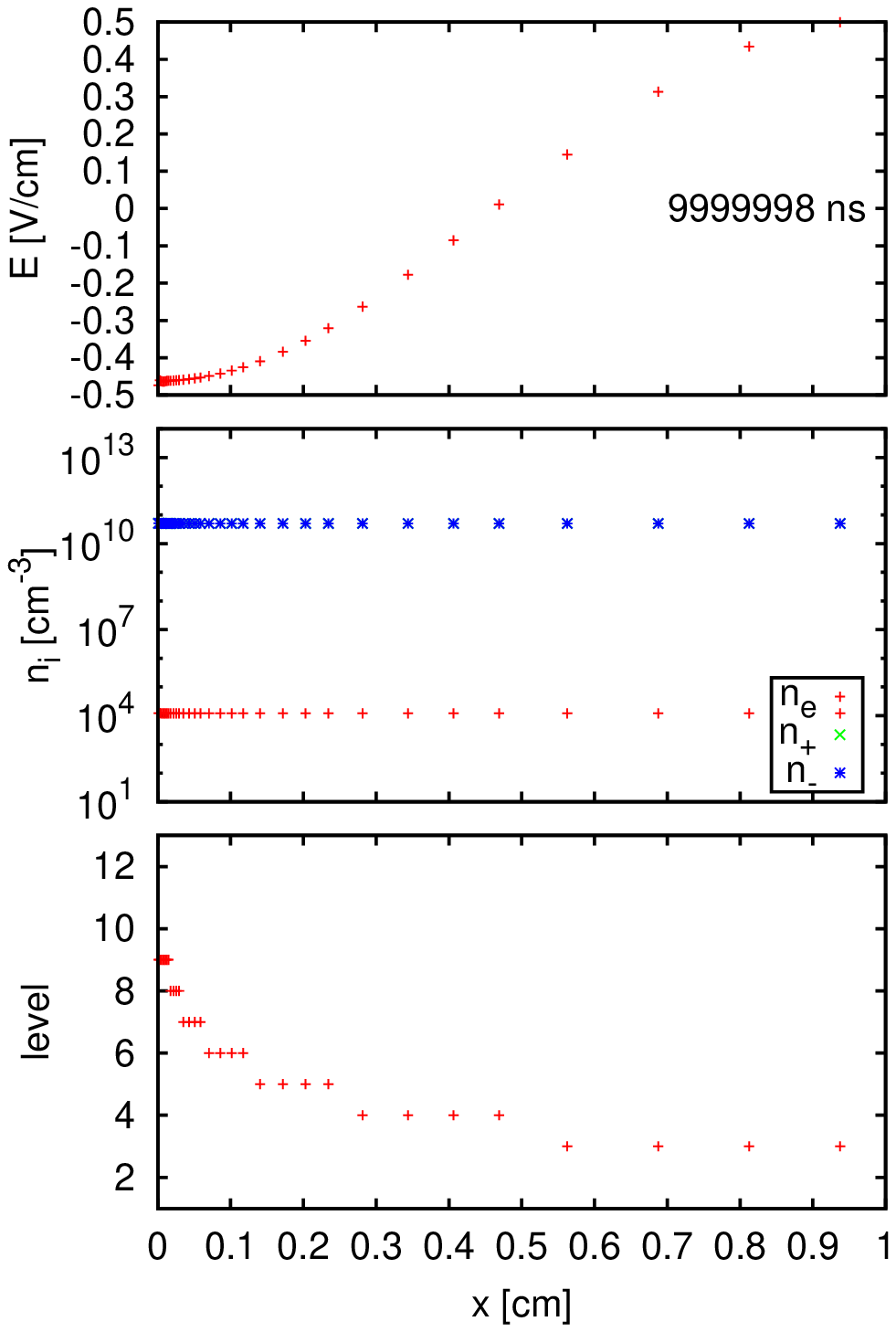}
   \caption{Steady-state of pulsed discharges (last period).
           Top: propagation of the discharge in the domain at $t=9900010\,$ns after the beginning of the pulse (left);
            and at $t=9900015\,$ns (right). 
            Bottom: relaxation on short time scale $t=9900050\,$ns;
            and end of the relaxation phase $t=9999998\,$ns (right).
 \label{multipulse-domain-3}}
 \end{center}
\end{figure}

After a few repetitive pulses, we have observed that the discharge dynamics reached a steady-state behavior as observed 
in the experiments. To show the characteristics of the discharge when the steady-state is reached, 
Figure~\ref{multipulse-domain-3} shows the discharge dynamics of the 100th period. The sequence of images is the 
same as in Figure~\ref{multipulse-domain-1}. 
At the end of the 99th pulse, we have observed that the axial distribution of charged species in the gap is uniform 
and that the level of preionization is $5 \times 10^{10}\,$cm$^{-3}$ positive and negative ions and $10^{4}\,$cm$^{-3}$ electrons.
We note that $10\,$ns after the beginning of the 100th pulse the propagation of the discharge is faster than 
for the first pulse. This faster propagation is mostly due to the higher preionization level of positive and negative 
ions in the gap in comparison of the first voltage pulse. 
We observe that for the 100th pulse, 
$15\,$ns after the beginning of the pulse the discharge has almost completely
crossed the inter-electrode gap and then 
during the relaxation phase, there is no remaining space charge in the whole gap. 
Consequently, $50\,$ns after the beginning of the 100th pulse, axial distributions of all charged species are uniform. 
As already observed for the first pulse, at $50\,$ns after the beginning of the voltage pulse most electrons 
have attached and then, the density of positive ions is almost equal to the density of negative ions in the whole gap. 
We see that the corresponding electric field distribution is not uniform at $50\,$ns, but no steep gradients are 
observed as for the first voltage pulse. At $t=9999998\,$ns, that is to say at the end of the 100th period, we note 
that a very low electric field is obtained in the gap. An axially uniform distribution of charges is obtained 
with $5 \times 10^{10}\,$cm$^{-3}$ for positive and negative ions and $10^{4}\,$cm$^{-3}$ for electrons, 
which was the initial condition of the 100th pulse.
This demonstrates the existence of a steady-state 
behavior of these nanosecond repetitively pulsed discharges.

\FloatBarrier

\section{Conclusions}
The present work proposes a new numerical 
strategy for multi-scale streamer simulations.
It is based on an adaptive second order time
integration strategy that allows to discriminate
time scales-related features of the phenomena,
given a required level of accuracy of computations.
Compared with a standard procedure for which accuracy is guaranteed
by considering time steps of the order of the fastest scale,
the control error approach
implies on the one hand, 
an effective accurate resolution independent of the
fastest physical time scale,
and on the other hand, an important improvement of
computational efficiency whenever the required time steps
go beyond standard stability constraints.
The latter is a direct consequence of the self-adaptive time step
strategy for the resolution of the drift-diffusion equations
which considers
splitting time steps 
not limited by stability constraints for reaction, diffusion and convection
phenomena.
So far, the global decoupling time steps are limited by the dielectric
relaxation stability constraint but with a second order accuracy.
Nevertheless, we have also demonstrated that
the decoupling time steps are rather chosen based on an accuracy criterion.
Besides,
if a technique such as a semi-implicit approach is implemented,
the same ideas of the proposed adaptive strategy remain valid.

An adaptive multiresolution
technique was also proposed in order to provide error control of
the spatial adapted representation.
The numerical results have proven a natural coupling between time
and space accuracy requirements and how the set of time-space accuracy
tolerances tunes the precise description of the overall time-space
multi-scale phenomenon.
As a consequence, the numerical results for
multi-pulsed discharge configurations prove 
that
this kind of multi-scale phenomena, previously out of reach, can be 
successfully simulated with conventional computing resources
by this time-space adaptive strategy.
And they also show that 
a consistent physical description
is achieved
for a broad spectrum of space and time scales 
as well as different physical scenarios.

In this work, we focused on a 1.5D model
in order to
evaluate the numerical performance of the strategy.
However, 
the dimension of the problem will only have an influence
on the computational efficiency measurements but not on any space-time accuracy 
or stability aspects.
At this stage of development, the same numerical strategy can be
coupled with a multi-dimensional Poisson's equation solver,
even for adapted grid configurations as developed 
recently in \cite{Montijn:2006,Unfer:2010,Pancheshnyi:2008}.
Finally,
an important amount of work is still in progress
concerning 
programming features such as data structures,
optimized routines and parallelization
strategies. 
On the other hand,
numerical analysis
of theoretical aspects is also underway
to extend and further improve the proposed numerical strategy.
These issues constitute particular topics of our current research.

%\FloatBarrier

%\section*{References}

\bibliographystyle{elsarticle-num}
\bibliography{plasma1d.bib}

\end{document}